\documentstyle[11pt]{article}
\title{The connectedness of some varieties and the Deligne-Simpson problem}
\author{Vladimir Petrov Kostov\\ \\ \hspace{7cm}
{\sl To the memory of my mother}\\ \\ } 
\date{}
\newtheorem{tm}{Theorem}
\newtheorem{lm}[tm]{Lemma}
\newtheorem{cor}[tm]{Corollary}
\newtheorem{prop}[tm]{Proposition}
\newtheorem{rem}[tm]{Remark}

\newtheorem{defi}[tm]{Definition}
\newtheorem{ex}[tm]{Example}
\newtheorem{conjecture}[tm]{Open questions and comments}
\newtheorem{nota}[tm]{Notation}
\begin{document}
\maketitle 

\begin{abstract}
The Deligne-Simpson problem (DSP) (resp. the weak DSP)   
is formulated like this: {\em give necessary and 
sufficient conditions for the choice of the conjugacy classes 
$C_j\subset GL(n,{\bf C})$ or $c_j\subset gl(n,{\bf C})$ so that there exist 
irreducible (resp. with trivial centralizer) $(p+1)$-tuples  
of matrices $M_j\in C_j$ or $A_j\in c_j$ 
satisfying the equality $M_1\ldots M_{p+1}=I$ or $A_1+\ldots +A_{p+1}=0$}. 
The matrices $M_j$ and $A_j$ are interpreted as monodromy operators of 
regular linear systems and as matrices-residua of Fuchsian ones on Riemann's 
sphere. For $(p+1)$-tuples of conjugacy classes one of 
which is with distinct eigenvalues 1) we prove that the variety 
$\{ (M_1,\ldots ,M_{p+1})|M_j\in C_j,M_1\ldots M_{p+1}=I\}$ or 
$\{ (A_1,\ldots ,A_{p+1})|A_j\in c_j,A_1+\ldots +A_{p+1}=0\}$ is connected if 
the DSP is positively solved for the given conjugacy classes and 2)  
we give necessary and sufficient conditions for the positive solvability 
of the weak DSP.\\  


{\bf AMS classification index:} 15A30, 15A24, 20G05
\end{abstract}

\section{Introduction}

The {\em Deligne-Simpson problem (DSP)} is formulated like this: {\em give 
necessary and sufficient conditions upon the choice of the conjugacy classes 
$c_j\subset gl(n,{\bf C})$ (resp. $C_j\subset GL(n,{\bf C})$), 
$j=1,\ldots ,p+1,$, $p\geq 2$, 
so that there exist irreducible $(p+1)$-tuples of matrices 
$A_j\in c_j$ (resp. $M_j\in C_j$) such that $A_1+\ldots +A_{p+1}=0$ 
(resp. $M_1\ldots M_{p+1}=I$).} 
``Irreducible'' means ``without common proper invariant subspace''. 
The {\em weak DSP} is 
obtained from the DSP by replacing the requirement of 
irreducibility by the weaker requirement the centralizer of the 
$(p+1)$-tuple to be trivial, i.e. reduced to scalars. In what follows we 
write ``tuple'' instead of ``$(p+1)$-tuple''.

The DSP or the weak DSP is 
{\em solvable} for given conjugacy classes $c_j$ or $C_j$ if there exist 
tuples of matrices $A_j\in c_j$ whose sum is $0$, resp. 
of matrices $M_j\in C_j$ whose product is $I$, irreducible or with trivial 
centralizer. 

The matrices $A_j$ (resp. $M_j$) are interpreted as {\em matrices-residua} 
(resp. as {\em monodromy matrices}) of Fuchsian (resp. regular) linear 
systems of differential equations on 
Riemann's sphere; the conjugacy classes $C_j$ are interpreted as local 
monodromies around the poles of a regular linear system and 
for matrices $M_j$ the problem admits the interpretation: {\em for which 
tuples of local monodromies do there exist irreducible monodromy 
groups with such local monodromies}; 
see the details in \cite{Ko1} or \cite{Ko2}. 

The multiplicative version of the problem (i.e. for matrices $M_j$) has been 
stated by P. Deligne (the additive, i.e. for matrices $A_j$, by the author) 
and C. Simpson was the first to obtain significant results towards the 
resolution of the problem, see \cite{Si}.

In what follows we assume that there 
hold the self-evident necessary conditions the sum of the traces of the 
classes $c_j$ to be $0$ (resp. the product of the determinants of the 
classes $C_j$ to be $1$). In terms of the eigenvalues $\lambda _{k,j}$  
(resp. $\sigma _{k,j}$) of the matrices from $c_j$ (resp. $C_j$) repeated 
with their multiplicities, this condition reads     

\[ \sum _{j=1}^{p+1}\sum _{k=1}^n\lambda _{k,j}=0~~  
{\rm resp.}~~\prod _{j=1}^{p+1}\prod _{k=1}^n\sigma _{k,j}=1~.\] 

\begin{defi}
An equality  
$\sum _{j=1}^{p+1}\sum _{k\in \Phi _j}\lambda _{k,j}=0$, resp. 
$\prod _{j=1}^{p+1}\prod _{k\in \Phi _j}\sigma _{k,j}=1$, is called a 
{\em non-genericity relation};  
the non-empty sets $\Phi _j$ contain one and the same number $<n$ of indices  
for all $j$. Eigenvalues satisfying none of these relations are called 
{\em generic}.\end{defi} 

\begin{rem}
Reducible  
tuples exist only for non-generic eigenvalues (the 
eigenvalues of each diagonal block of a block upper-triangular  
tuple satisfy some non-genericity relation).
\end{rem}   

For generic eigenvalues the 
problem is completely solved in \cite{Ko1}, \cite{Ko2} and \cite{Ko3}, 
and the solution is a 
criterium upon the Jordan normal forms (JNFs) defined by the conjugacy 
classes, i.e. 
it does not depend on the concrete choice of the eigenvalues provided that 
the latter remain generic. The additive version of the problem is solved also 
for arbitrary eigenvalues, see \cite{C-B}.

Denote by $J(Y)$ (resp. $J(C)$) the JNF of the matrix $Y$ (resp. the JNF 
defined by the conjugacy class $C$). For a conjugacy class $C$ in 
$gl(n,{\bf C})$ or $GL(n,{\bf C})$ denote by 
$d(C)$ its dimension (as a variety in $gl(n,{\bf C})$ or $GL(n,{\bf C})$) 
and for a matrix $Y\in C$ set 
$r(C):=\min _{\lambda \in {\bf C}}{\rm rank}(Y-\lambda I)$. The integer 
$n-r(C)$ is the maximal number of Jordan blocks of $J(Y)$ with one and the 
same eigenvalue. Set $d_j:=d(c_j)$ (resp. $d(C_j)$), $r_j:=r(c_j)$ 
(resp. $r(C_j)$). The quantities 
$r(C)$ and $d(C)$ depend only on the JNF $J(C)$, not 
on the eigenvalues. 

\begin{prop}\label{d_jr_j}
(C. Simpson, see \cite{Si}.) The 
following couple of inequalities is a necessary condition for solvability 
of the DSP in the case of matrices $M_j$:

\[ d_1+\ldots +d_{p+1}\geq 2n^2-2~~~~~(\alpha _n)~~,~~~~~
{\rm for~all~}j,~r_1+\ldots +\hat{r}_j+\ldots +r_{p+1}\geq n~~~~~
(\beta _n)~~~.\]
\end{prop}

The proposition holds also for matrices $A_j$ whose sum is $0$, 
see \cite{Ko2} or \cite{Ko4}. 

\begin{prop}\label{trivcentr}
Conditions $(\alpha _n)$ and $(\beta _n)$ are necessary for the solvability 
of the weak DSP in the case when one of the conjugacy classes has $n$ distinct 
eigenvalues.
\end{prop}

The proposition is proved at the end of part B) of Subsection~\ref{knownfacts}.
C. Simpson proves in \cite{Si} the following 

\begin{tm}\label{CSimpson}
For generic eigenvalues and when one of 
the conjugacy classes $C_j$ is with distinct eigenvalues, conditions 
$(\alpha _n)$ and $(\beta _n)$ together are necessary and sufficient for the 
existence of irreducible tuples of matrices $M_j\in C_j$ whose product 
is $I$.
\end{tm}

The same is true for matrices $A_j\in c_j$ whose sum is $0$ (see 
\cite{Ko4} and compare with \cite{Ko5}, Theorems~19 and 32 in which one of 
the matrices is supposed to have eigenvalues of multiplicity $\leq 2$, not 
necessarily distinct ones).  

In the present paper we consider sets of the form 

\[ {\cal V}(c_1,\ldots c_{p+1})=
\{ (A_1,\ldots ,A_{p+1})|A_j\in c_j, A_1+\ldots +A_{p+1}=0\} \]

and    
\[ {\cal W}(C_1,\ldots C_{p+1})=
\{ (M_1,\ldots ,M_{p+1})|M_j\in C_j, M_1\ldots M_{p+1}=I\} \]
or just ${\cal V}$ and ${\cal W}$ for short. 
The aim of the present paper is to prove (in Section~\ref{prtm}) the following 

\begin{tm}\label{basictm}
1) For generic eigenvalues and when one of 
the conjugacy classes $c_j$ (resp. $C_j$) is with distinct eigenvalues the 
set ${\cal V}$ (resp. ${\cal W}$) 
is a smooth and connected variety.

2) For arbitrary eigenvalues, when one of 
the conjugacy classes $c_j$ (resp. $C_j$) is with distinct eigenvalues, and if 
there exist irreducible tuples, then the 
closure of the set ${\cal V}$ (resp. 
${\cal W}$) 
is a connected variety. The algebraic closures of these sets coincide with 
their topological closures; these are the closures  
of the subvarieties consisting of irreducible 
tuples; these subvarieties are connected. 
The singular points of the closures are precisely the tuples 
of matrices with non-trivial centralizers. 

3) If one of the matrices $A_j$ or $M_j$ is with distinct eigenvalues, then 
conditions $(\alpha _n)$ and $(\beta _n)$ together are necessary and 
sufficient for the solvability of the weak DSP.
\end{tm}

{\bf Convention.} {\em In what follows we assume that the conjugacy class 
$c_{p+1}$ (resp. $C_{p+1}$) is with distinct eigenvalues.}

\begin{ex}\label{threestrata}
Consider the 
case of matrices $A_j$ for $p=2$, 
$n=2$, with eigenvalues of the three diagonalizable non-scalar 
matrices respectively $(a,b)$, $(c,d)$, 
$(g,h)$, where $a+c+g=b+d+h=0$ and there are no non-genericity 
relations which are not corollaries of these ones. Then the set ${\cal V}$ 
is a stratified variety with three strata -- $S_0$, $S_1$ and $S_2$. The 
stratum $S_i$ consists of triples which up to conjugacy equal 

\[ A_1=\left( \begin{array}{cc}a&0\\0&b\end{array}\right) ~~,~~
A_2=\left( \begin{array}{cc}c&\varepsilon _i\\ \eta _i&d
\end{array}\right) ~~, ~~
A_3=\left( \begin{array}{cc}g&-\varepsilon _i\\-\eta _i&h\end{array}\right) \]
where $\varepsilon _0=\eta _0=0$, $\varepsilon _1=1$, $\eta _1=0$, 
$\varepsilon _2=0$, $\eta _2=1$. Hence, $S_0$ lies in the closure of $S_1$ 
and $S_2$. A theorem by N. Katz (see \cite{Ka}) forbids coexistence of 
irreducible and reducible triples in the so-called {\em rigid} case (i.e. 
when there is equality in $(\alpha _n)$ which is the case here -- one has 
$d_1=d_2=d_3=2$). Hence, there 
are no strata of ${\cal V}$ other than $S_0$, $S_1$ and $S_2$, and ${\cal V}$ 
is connected. However, it is not smooth at $S_0$. Indeed, one can deform 
analytically a triple from $S_0$ into one from $S_1$ and into one from $S_2$; 
the triples from $S_1$ and $S_2$ defining different semi-direct sums, the 
strata $S_1$ and $S_2$ cannot be parts of one and the same smooth variety 
containing $S_0$.  
\end{ex}      

\begin{ex}
Again for $p=n=2$, consider the case when $c_1$ is nilpotent non-scalar and 
$c_2=-c_3$ is with eigenvalues $1$, $2$. Then the algebraic and topological 
closure of the variety 
${\cal V}(c_1,c_2,c_3)$ consists of three strata -- $T_0$, $T_1$ and $T_2$. Up 
to conjugacy, the triples from $T_0$ are diagonal and equal 
diag$(0,0)$, diag$(1,2)$, diag$(-1,-2)$. The ones from $T_1$ (resp. $T_2$) 
up to conjugacy equal  

\[ A_1=\left( \begin{array}{cc}0&1\\0&0\end{array}\right) ~~,~~
A_2=\left( \begin{array}{cc}1&-1\\0&2
\end{array}\right)  ~~({\rm resp.}~A_2=
\left( \begin{array}{cc}2&-1\\0&1\end{array}\right) )~~,~~A_3=-A_1-A_2.\]
The stratum $T_0$ lies in the closures of $T_1$ and $T_2$ and like in 
Example~\ref{threestrata}, the closure of 
${\cal V}$ is singular along $T_0$. Notice that 
${\cal V}$ itself is not connected because 
$T_0\not\in {\cal V}$ and like in Example~\ref{threestrata}, the strata 
$T_1$, $T_2$ are not parts of one and the same smooth variety.
\end{ex}

The reader will find other examples illustrating the stratified structure of 
the varieties ${\cal V}$ or 
${\cal W}$ in \cite{Ko7} and \cite{Ko8}, in particular, 
cases when the dimension of the variety is higher than the expected one when  
the centralizer is non-trivial.

\begin{conjecture}
1) It would be natural to ask the question whether parts 1) and 2) of  
Theorem~\ref{basictm} are true without the condition one of the classes $c_j$ 
or $C_j$ to be with distinct eigenvalues. The author is convinced that this 
is true (in the case when inequality $(\alpha _n)$ becomes equality the proof 
of this can be deduced from the results in \cite{Ka} -- rigid tuples are 
unique up to conjugacy and coexistence of irreducible and reducible tuples is 
impossible in the rigid case). 
In the present paper we use the condition one of the classes to be with 
distinct eigenvalues 
in the proofs (see Proposition~\ref{directsum} and its proof).

2) It would be interesting to prove the connectedness of the closures 
of the varieties 
${\cal V}$ or ${\cal W}$ without the assumption that there are irreducible 
tuples, see part 2) of Theorem~\ref{basictm}. All examples known to 
the author are of connected closures, see \cite{Ko7} and \cite{Ko8}.

3) In part 3) of the theorem the condition one of the matrices to be with 
distinct eigenvalues is essential. Even for double eigenvalues there is a 
counterexample -- for $n=2$ a triple of nilpotent non-zero matrices $A_j$ 
whose sum is $0$, is upper-triangular up to conjugacy, and, hence, the 
centralizer of the triple contains all three matrices. 
Hence, for $n=2$ the weak DSP is not solvable for 
a triple of nilpotent non-zero conjugacy classes.

4) The theorem implies that the moduli space of tuples of matrices 
from given conjugacy classes (one of which is with distinct eigenvalues),  
with zero sum or whose product is $I$, is connected provided that there exist 
irreducible tuples.
\end{conjecture}

\section{Preparation for the proof of Theorem \protect\ref{basictm}
\protect\label{preparation}}
\subsection{The known facts we use\protect\label{knownfacts}}

{\bf A) The dimension of a variety ${\cal V}$ or 
${\cal W}$.} 

To prove the theorem we need the following propositions: 

\begin{prop}\label{[A_j,X_j]}
The centralizer of the $p$-tuple of matrices $A_j$ ($j=1,\ldots ,p$) is 
trivial if and only if the mapping 
$(sl(n,{\bf C}))^p\rightarrow sl(n,{\bf C})$, 
$(X_1,\ldots ,X_p)\mapsto \sum _{j=1}^p[A_j,X_j]$ is surjective.
\end{prop}

The proposition is proved in \cite{Ko2}.

\begin{prop}\label{smooth}
At a point defining a tuple with trivial centralizer a variety 
${\cal V}$ or ${\cal W}$ is smooth and locally of dimension 
$d_1+\ldots +d_{p+1}-n^2+1$.
\end{prop}

The propositions of this subsection (except Propositions~\ref{[A_j,X_j]} 
and \ref{samedr}) are proved in the next one.\\   

{\bf B) The basic technical tool.}

\begin{defi}
Call {\em basic technical tool} the procedure described 
below whose aim is to deform analytically a given tuple of matrices 
$A_j$ or $M_j$ with trivial centralizer 
by changing their conjugacy classes in a desired way. 
\end{defi}

Set $A_j=Q_j^{-1}G_jQ_j$, 
$G_j$ being Jordan matrices. 
Look for a tuple of matrices $\tilde{A}_j$ (whose sum is $0$) of the form 

\[ \tilde{A}_j=(I+\varepsilon X_j(\varepsilon ))^{-1}
Q_j^{-1}
(G_j+\varepsilon V_j(\varepsilon ))Q_j
(I+\varepsilon X_j(\varepsilon ))\] 
where $\varepsilon \in ({\bf C},0)$ and 
$V_j(\varepsilon )$ are given 
matrices analytic 
in $\varepsilon$; they must satisfy the condition  
tr$(\sum _{j=1}^{p+1}V_j(\varepsilon ))
\equiv 0$; set $N_j=Q_j^{-1}V_jQ_j$. The existence of matrices 
$X_j(\varepsilon )$ is deduced from 
the triviality of the centralizer, using Proposition~\ref{[A_j,X_j]}  
(see its proof and the details in \cite{Ko2}).

Notice that one has 
$\tilde{A}_j=A_j+\varepsilon [A_j,X_j(0)]+\varepsilon N_j+o(\varepsilon )$. 
Proposition~\ref{[A_j,X_j]} assures the existence of $X_j(0)$, i.e. 
the existence in first approximation w.r.t. 
$\varepsilon$ of the matrices $X_j$; the existence of true matrices $X_j$ 
analytic in $\varepsilon$ follows from the implicit function theorem. 

If for $\varepsilon \neq 0$ 
small enough the eigenvalues of the matrices $\tilde{A}_j$ are generic, then 
their tuple is irreducible. In a similar way one can deform 
analytically tuples depending on a multi-dimensional parameter.

Deforming analytically tuples of matrices $A_j$ is used sometimes to make 
some of the JNFs $J(A_j)$ ``more generic'':

\begin{ex}\label{JNFnongeneric}
Suppose that $A_{p+1}$ is diagonal and 
has $n-2$ simple eigenvalues and a double eigenvalue (the tuple being with 
trivial centralizer). Set 

\[ A_{p+1}=\left( \begin{array}{ccc}a&0&0\\0&a&0\\0&0&D
\end{array}\right) ~~,~~N_{p+1}=\left( \begin{array}{ccc}0&1&0\\0&0&0\\0&0&0
\end{array}\right) \]
where $D$ is $n\times n$, diagonal, with distinct and different 
from $a$ eigenvalues. Then for $\varepsilon \neq 0$ the matrix $A_{p+1}$ is 
with the same eigenvalues and with a single Jordan block of size $2$ 
corresponding to the eigenvalue $a$.
\end{ex}

In other cases deforming analytically of tuples changes the eigenvalues 
of the matrices $A_j$:

\begin{ex}\label{newEVs}
Suppose that 

\[ A_{p+1}=\left( \begin{array}{ccc}a&1&0\\0&a&0\\0&0&D
\end{array}\right) ~~,~~N_{p+1}=\left( \begin{array}{ccc}0&0&0\\
f&g&0\\0&0&0
\end{array}\right) \]
where $D$ is as in Example~\ref{JNFnongeneric} and $f,g\in {\bf C}$. For 
generic $f$ and $g$ one obtains for $\varepsilon \neq 0$ small enough 
matrices $\tilde{A}_{p+1}$ with $n$ distinct eigenvalues. 
\end{ex}

Given a tuple of matrices $M_j$ with  
trivial centralizer and whose product is $I$, look for matrices 
$\tilde{M}_j$ (whose product is $I$) of the form 

\[ \tilde{M}_j=(I+\varepsilon X_j(\varepsilon ))^{-1}(M_j+
\varepsilon N_j(\varepsilon ))(I+
\varepsilon X_j(\varepsilon ))\] 
where the given matrices $N_j$ depend analytically on 
$\varepsilon \in ({\bf C},0)$ and the product of the determinants of the 
matrices $\tilde{M}_j$ is $1$; one looks for $X_j$ analytic in 
$\varepsilon$. The existence of 
such matrices $X_j$ follows again from the triviality of the centralizer, 
see \cite{Ko2}.

When applying the basic technical tool one often preserves the conjugacy 
classes of all matrices $A_j$ or $M_j$ but one, as this is done in the 

{\em Proof of Proposition~\ref{trivcentr}:}

Given a tuple with trivial centralizer 
of matrices $A_j$ whose sum is $0$ or of matrices 
$M_j$ whose product is $I$ the matrix $A_{p+1}$ or $M_{p+1}$ having distinct 
eigenvalues, deform it into a nearby such tuple with generic 
eigenvalues, hence, irreducible; for $j=1,\ldots ,p$ the matrix $A_j$ 
or $M_j$ remains within its conjugacy class; $A_{p+1}$ (resp. $M_{p+1}$) 
remains with distinct eigenvalues, but the latter change. For the deformed 
tuple conditions $(\alpha _n)$ and $(\beta _n)$ hold, hence, 
they hold for the initial 
one as well (the quantities $d_j$, $r_j$ do not change under such a 
deformation and when checking condition $(\beta _n)$ one needs to 
consider only the sum $r_1+\ldots +r_p$ because $r_{p+1}=n-1$).~~~~$\Box$     
 
{\bf C) Corresponding Jordan normal forms.} 

A JNF is a collection of numbers $\{ b_{i,l}\}$ where $b_{i,l}$ 
stands for the size of the $i$-th Jordan block with the $l$-th 
eigenvalue. For a given JNF  
$J^n=\{ b_{i,l}\}$ of size $n$ define its {\em corresponding} 
diagonal JNF ${J'}^n$ (more generally, two JNFs are said to be 
{\em corresponding} to one another if they correspond to one and the same 
diagonal JNF). A diagonal JNF is  
a partition of $n$ defined by the multiplicities of the eigenvalues. 
For each $l$ fixed, the collection of numbers 
$\{ b_{i,l}\}$ is a partition of $\sum _{i\in I_l}b_{i,l}$ and 
${J'}^n$ is the disjoint sum of the dual partitions; 
see the details in \cite{Ko1}, \cite{Ko2} or \cite{Kr}, as well as the proof 
of the following 

\begin{prop}\label{samedr}
The quantities $d$ and $r$ 
defined in the Introduction are the same for two corresponding JNFs.
\end{prop}
 
{\bf D) Realizing of monodromy groups by Fuchsian systems with different sets 
of eigenvalues of their matrices-residua. Procedure $(l,k)$.} 

\begin{nota}\label{matrixentry}
$E_{i,k}$ denotes the 
matrix with a single non-zero entry which equals $1$ and is 
in position $(i,k)$. In all other cases double subscripts denote matrix 
entries. E.g., $A_{j;1,2}$ stands for the 
entry of the matrix $A_j$ in position $(1,2)$. 
\end{nota}

When one makes the linear change of variables $X\mapsto V(t)X$ (where $V(t)$ 
is an $n\times n$-matrix-function meromorphic on the Riemann sphere and whose 
determinant does not vanish identically there) in the linear 
system d$X/$d$t=A(t)X$, then the matrix $A(t)$ undergoes the gauge 
transformation: 
$A\mapsto -V^{-1}$d$V/$d$t+V^{-1}AV$. 
If the system is Fuchsian, then after such a 
change, in general, it is no longer Fuchsian, but for special 
choices of $V$ one can obtain 
new Fuchsian systems (with the same monodromy group) whose matrices-residua 
belong to new conjugacy classes. For each matrix-residuum the eigenvalues 
of the new conjugacy class are shifted by integers w.r.t. the ones of the 
old conjugacy class. Indeed, the eigenvalues $\lambda _{k,j}$ of $A_j$ and 
$\sigma _{k,j}$ of $M_j$ are related by the equality  
$\exp (2\pi i\lambda _{k,j})=\sigma _{k,j}$. 

We are interested here only in changes that shift the eigenvalues of $A_{p+1}$ 
while preserving the conjugacy classes of the other matrices $A_j$. This is 
the case when $V$ is holomorphic and holomorphically invertible for 
$t\neq a_{p+1}$. An {\em admissible shift} of the eigenvalues of $A_{p+1}$ 
is a shift by an $n$-vector with integer components whose sum is $0$. 
The basic result from \cite{Ko6} implies the following 

\begin{cor}\label{QS}
If the monodromy group defined by the operators $M_j$ is irreducible, then 
for all but finitely many admissible shifts there exist linear changes 
$X\mapsto V(t)X$ transforming the given Fuchsian system into a new one 
with the eigenvalues of $A_{p+1}$ shifted like required.
\end{cor} 

Explain how to perform the simplest admissible shifts. 
Consider the Fuchsian system 

\begin{equation}\label{Fuchs}
{\rm d}X/{\rm d}t=(\sum _{j=1}^{p+1}A_j/(t-a_j))X
\end{equation}
Its Laurent series expansion at $a_{p+1}$ looks like this: 

\begin{equation}\label{Fuchs1}
{\rm d}X/{\rm d}t=(A_{p+1}/(t-a_{p+1})+\sum _{j=1}^pA_j/(a_{p+1}-a_j)+
o(1))X
\end{equation}
Assume that $A_{p+1}={\rm diag}(\lambda _{1,p+1},\ldots ,\lambda _{n,p+1})$. 
If the entry $c=(\sum _{j=1}^pA_j/(a_{p+1}-a_j))_{l,k}$ is non-zero, then in 
the Fuchsian system (\ref{Fuchs}) one 
can perform the change $w:X\mapsto (I+W/(t-a_{p+1}))X$ (with 
$W_{k,l}=(\lambda _{k,p+1}-\lambda _{l,p+1}+1)/c$, $l\neq k$, $W_{i,j}=0$ for 
$(i,j)\neq (k,l)$) 
after which the eigenvalues of $A_{p+1}$ change as follows: 

\[ \lambda _{l,p+1}\mapsto \lambda _{l,p+1}-1~,~
\lambda _{k,p+1}\mapsto \lambda _{k,p+1}+1~,~
\lambda _{i,p+1}\mapsto \lambda _{i,p+1}~{\rm for~}i\neq l,k~.\]
Observe that the matrix $I+W/(t-a_{p+1})$ is holomorphic and holomorphically 
invertible for $t\neq a_{p+1}$ and that its determinant is identically 
equal to $1$. Call this change of eigenvalues {\em Procedure $(l,k)$}.

\begin{prop}\label{Procedure}
Given a system (\ref{Fuchs}) with irreducible tuple of 
matrices-residua $A_j\in c_j$, 
$A_{p+1}={\rm diag}(\lambda _{1,p+1},\ldots ,\lambda _{n,p+1})$, and given 
$l\neq k$, $1\leq l,k\leq n$, one can either perform Procedure $(l,k)$ or 
one can change the matrices $A_j$, $1\leq j\leq p$, within their 
conjugacy classes (their sum remaining the same) and/or the 
positions of the poles $a_j$ so that Procedure $(l,k)$ will be possible 
to perform in the changed system. The new tuple of poles $a_j$ 
and the new $p$-tuple of matrices $A_j$, $1\leq j\leq p$, can be chosen 
arbitrarily close to the initial ones.
\end{prop}

{\bf E) The mapping $\phi$.}

For a matrix $A_j\in c_j$ set $M(A_j)=\exp (2\pi iA_j)$. 
Denote by $C_j$ the conjugacy class of the matrix  $M(A_j)$. The 
monodromy operators $M_j$ of the Fuchsian system (\ref{Fuchs}) defined 
after a standard set of 
generators of $\pi _1({\bf C}P^1\backslash \{ a_1,\ldots ,a_{p+1}\} )$ (see 
the details in \cite{Ko2}) equal $M(A_j)$ provided that there is no non-zero 
integer difference between two eigenvalues of $A_j$. We assume that a fixed 
point $a_0$ is chosen as well as an initial value of the solution, i.e. a 
non-degenerate $n\times n$-matrix $B=X|_{t=a_0}$. Thus the mapping 
$\phi :(A_1,\ldots ,A_{p+1})\mapsto (M_1,\ldots ,M_{p+1})$ is a mapping 
${\cal V}(c_1,\ldots ,c_{p+1})\rightarrow {\cal W}(C_1,\ldots ,C_{p+1})$.      

\begin{prop}\label{VmapstoW}
1) Suppose that no conjugacy class $c_j$ has two eigenvalues differing by a 
non-zero integer. Then the mapping $\phi$ is a local diffeomorphism at every 
point of ${\cal V}$ where the centralizer of the 
tuple of matrices $A_j$ is trivial. 

2) If in addition the eigenvalues of the conjugacy classes $c_j$ are generic, 
then $\phi$ is a global 
diffeomorphism of ${\cal V}$ onto its image in 
${\cal W}$.
\end{prop}

{\bf F) Deforming direct sums into semi-direct ones and into irreducible 
representations.}

\begin{prop}\label{directsum}
1) A tuple of matrices $A_j$ from the closure of 
${\cal V}$ (or of ${\cal W}$) 
and with non-trivial centralizer defines a 
representation which is a direct sum.

2) If $A_{p+1}$ (or $M_{p+1}$) is with distinct eigenvalues, 
then any matrix from the centralizer of a tuple of the variety 
${\cal V}$ (or ${\cal W}$) is 
diagonalizable.
\end{prop}

Denote by $P_i$, $i=1,2$, representations with trivial centralizers defined 
by matrices $A_j^i$ whose sum is $0$ or by matrices $M_j^i$ whose product is 
$I$; $j=1,\ldots ,p+1$. The two matrices $A_{p+1}^i$ (resp. $M_{p+1}^i$), 
$i=1,2$, are each with distinct eigenvalues and have no eigenvalue in 
common. Denote by $d_j^i$ (resp. $d_j'$) the quantities $d$ of the 
conjugacy classes of the matrices $A_j^i$ or $M_j^i$ (resp. of the direct 
sums $A_j^1\oplus A_j^2$ or $M_j^1\oplus M_j^2$). 
Set $\xi =$dim~Ext$^1(P_1,P_2)$.

\begin{prop}\label{dirsemidir}
1) One has $\xi =(\sum _{j=1}^{p+1}(d_j'-d_j^1-d_j^2)/2)-2m_1m_2$.

2) If $\xi \geq 1$, then there exists a semidirect sum of $P_1$ 
and $P_2$ which is not reduced to a direct one.

3) If the representations $P_i$ are irreducible and $\xi \geq 2$, then the 
semidirect sum of $P_1$ and $P_2$ can be deformed into an irreducible 
representation without changing the conjugacy classes of the matrices 
$A_j$ or $M_j$ (their sum remaining $0$ or their product remaining $I$). 
\end{prop}

\begin{rem}\label{aposteriori}
One knows a posteriori (i.e. after Theorem~\ref{basictm} is proved) that 
in the conditions of the proposition, when $P_1$, $P_2$ are irreducible,  
if $\xi =1$ (and unlike for $\xi \geq 2$), then there exist no irreducible 
tuples. Indeed, in this case dim${\cal V}$ equals 
the dimension of its subvariety ${\cal U}$ of semi-direct sums (this can be 
proved like Proposition~\ref{dirsemidir}), so ${\cal U}$ cannot belong to 
the closure of the set of irreducible representations. The existence 
of irreducible representations would contradict 
part 2) of Theorem~\ref{basictm}.
\end{rem}
 
\subsection{Proofs of the propositions}

{\bf Proof of Proposition~\ref{smooth}:}

$1^0$. Consider first the case of matrices $A_j$. 
Assume without restriction that $c_j\subset sl(n,{\bf C})$.  
Consider the cartesian product 
$(c_1\times \ldots \times c_p)\subset (sl(n,{\bf C}))^p$.  
Define the mapping 
$\tau :(c_1\times \ldots \times c_p)\rightarrow sl(n,{\bf C})$ by 
the rule $\tau :(A_1,\ldots ,A_p)\mapsto A_{p+1}=-A_1-\ldots -A_p$ 
(recall that the sum of the matrices $A_j$ is $0$). 

$2^0$. The variety ${\cal V}$  
is the intersection of the two 
varieties in $c_1\times \ldots \times c_p\times sl(n,{\bf C})$: 
the graph of the mapping $\tau$ 
and $c_1\times \ldots \times c_p\times c_{p+1}$. This intersection is 
transversal which implies the smoothness of ${\cal V}$. 
Transversality follows from Proposition~\ref{[A_j,X_j]} -- 
the tangent space to the conjugacy class $c_j$ at $A_j$ 
equals $\{ [A_j,X]|X\in gl(n,{\bf C})\}$. 

$3^0$. Recall that dim$c_j$ is denoted by $d_j$. One has   

\[ {\rm dim}\, {\cal V}=
(\sum _{j=1}^pd_j)-{\rm codim}_{sl(n,{\bf C})}c_{p+1}=
(\sum _{j=1}^pd_j)-[(n^2-1)-
d_{p+1}]=\sum _{j=1}^{p+1}d_j-n^2+1~.\]

$4^0$. In the case of matrices $M_j$ the only difference in the proof  
is that the mapping    
$(A_1,\ldots ,A_p)$ $\mapsto A_{p+1}=-A_1-\ldots -A_p$
from ${\bf 2^0}$ has to be replaced by the mapping 

\[(M_1,\ldots ,M_p)\mapsto M_{p+1}=(M_1\ldots M_p)^{-1}~.\]  
The reader will be able to restitute the missing technical details after 
examining the more detailed description of the basic technical tool given in 
\cite{Ko2}. 

The proposition is proved.~~~~~$\Box$

{\bf Proof of Proposition \ref{Procedure}:}

$1^0$. Assume for convenience that $(l,k)=(1,n)$. Procedure $(1,n)$ is defined 
only if $c=(\sum _{j=1}^pA_j/(a_{p+1}-a_j))_{1,n}\neq 0$. Hence, 
if at least one of the entries $A_{j;1,n}$, $1\leq j\leq p$, is non-zero, 
one can change the positions of the 
poles $a_j$ a little to obtain the condition $c\neq 0$. So assume that for 
$j=1,\ldots ,p$ one has $A_{j;1,n}=0$. (For $n=2$ this means that the 
tuple is reducible which is a contradiction; so in what follows we 
assume that $n>2$.) We show that there exist infinitesimal   
conjugations $A_j\mapsto (I+\varepsilon X_j)^{-1}A_j(I+\varepsilon X_j)$, 
$1\leq j\leq p$, 
$\varepsilon \in ({\bf C},0)$, the result of which is that for at least one 
$j$ one will have $A_{j;1,n}\neq 0$ for $\varepsilon \neq 0$ 
and $A_{p+1}$ does not change (in what follows we set $j=1$). By analogy with 
the basic technical tool (Subsection~\ref{knownfacts}, part B)) one can 
show that then there exist true (not only infinitesimal) conjugations. 

$2^0$. These conjugations must be such that 

\begin{equation}\label{FA}
\sum _{j=1}^p[A_j,X_j]=0
\end{equation}
(first approximation w.r.t. $\varepsilon$). If no such $p$-tuple of 
conjugations changes the matrices $A_j$ so that $A_{1;1,n}\neq 0$ in 
first approximation w.r.t. $\varepsilon$, then the following condition must 
be a corollary of (\ref{FA}): tr$(E_{n,1}[A_1,X_1])=0$. This means that the 
linear form tr$(E_{n,1}[A_1,X_1])=$tr$([E_{n,1},A_1]X_1)$ 
(the variables are the entries of $X_1$) 
must be representable as 
tr$(D\sum _{j=1}^p[A_j,X_j])=$tr$(\sum _{j=1}^p[D,A_j]X_j)$ for some matrix 
$D\in sl(n,{\bf C})$. So from now on we want to prove that such a matrix $D$ 
does not exist. 

$3^0$. Suppose it does. Hence, one has 

\begin{equation}\label{D}
[D,A_j]=0~{\rm for~}j=2,\ldots ,p~~,~~[D,A_1]=[E_{n,1},A_1]
\end{equation}
Summing up these equalities one deduces that 

\begin{equation}\label{DA}
-[D,A_{p+1}]=[E_{n,1},A_1]
\end{equation} 
The matrix $[E_{n,1},A_1]$ has non-zero entries only in the first 
column and in the last row. Hence, the possible non-zero entries of $D$ 
are in the first column, in the last row and on the diagonal (recall that 
$A_{p+1}$ is diagonal and with distinct eigenvalues).

$4^0$. Two cases are possible:

{\em Case A)}. $D$ has at least two distinct eigenvalues.

{\em Case B)}. $D$ is nilpotent.      
 
In case A) permute the second, $\ldots$, $(n-1)$-st eigenvalues of $A_{p+1}$ 
(by conjugation with a permutation matrix which is block-diagonal, with 
diagonal blocks of sizes $1$, $n-2$, $1$) so that the first $m$ eigenvalues 
of $D$ be equal (say, equal to $b$) and different from the 
$(m+1)$-st, $\ldots$, $(n-1)$-st one. 

If $m=1$, then equalities 
(\ref{D}) imply that the entries of $A_1,\ldots ,A_p$ in positions $(1,i)$, 
$2\leq i\leq n-1,$ equal $0$. As $A_{j;1,n}=0$ for $j=1,\ldots ,p$, the 
tuple of matrices $A_j$ is reducible -- a contradiction. 

If $m>1$ and if $b$ is different from the $n$-th 
eigenvalue of $D$ as well, then it follows from (\ref{D}) that the 
entries of $A_1$, $\ldots$, $A_p$ in the right upper corner $m\times (n-m)$ 
are zeros. Hence, the tuple of matrices $A_j$ is reducible -- a 
contradiction.

If $m>1$ and the $n$-th eigenvalue of $D$ equals $b$, then block-decompose the 
$n\times n$-matrices as follows: $\left( \begin{array}{ccc}B&M&Q\\
S&R&T\\U&V&b\end{array}\right)$ where $B$ is $m\times m$ and $b$ is 
$1\times 1$. Conjugate the matrices $A_j$ and $D$ 
by a matrix $I+F$ (where only the $S$-, $V$- and $U$-blocks of $F$ are 
non-zero) to annihilate the $S$- and $V$-blocks of $D$; note that 
$(I+F)^{-1}E_{n,1}(I+F)=E_{n,1}$.   
After the conjugation equalities (\ref{D}) imply that the 
$M$- and $T$-blocks of the matrices $A_1,\ldots ,A_p$ are $0$. But then the 
$S$-blocks of these matrices are also $0$. Hence, the 
tuple is reducible -- a contradiction again. So case A) is impossible.

$5^0$. Consider case B). The form of $D$ (see $3^0$) implies that all entries 
of $D$ on and above the diagonal are $0$ and either 
$D^2=0$ and rk$D\leq 2$ or $D^2\neq 0$, rk$D=2$ and $D^3=0$. In 
the second case one has $D^2=\alpha E_{n,1}$, $\alpha \neq 0$ and it 
follows from (\ref{D}) that $[A_j,E_{n,1}]=0$ for $j=2,\ldots ,p$; hence, 
$A_{j;1,i}=0$ for $j=2,\ldots ,p$, $i=2,\ldots ,n$ and as these equalities 
hold also for $j=p+1$ ($A_{p+1}$ is diagonal), the tuple of 
matrices $A_j$ is reducible. 

So suppose that $D^2=0$ and rk$D\leq 2$. Let first rk$D=2$ (this is possible 
only for $n>3$). One can conjugate 
the matrices $A_j$ and $D$ by a block-diagonal matrix $W$ with three 
diagonal blocks, of sizes $1$, $n-2$, $1$, so that the matrix $D$ has at 
most three non-zero entries, in positions $(n-1,1)$, $(n,2)$ and $(n,1)$. 
The first two of them must be non-zero because rk$D=2$. It 
follows from (\ref{D}) that the entries in positions $(i,k)$ 
of the matrices $A_j$, $j=1,\ldots ,p$, are $0$ for $i=1,2$; $k=3,\ldots ,n$ 
(the details are left for the reader). Hence, the tuple is reducible.  

If rk$D=1$, then by a similar conjugation $D$ is simplified to have non-zero 
entries only in one or both of positions $(n-1,1)$, $(n,1)$ or 
$(n,1)$, $(n,2)$. In the first case the matrices $A_j$, $j=1,\ldots ,p$ have 
zero entries in positions $(1,2)$, $\ldots$, $(1,n)$, in the second case 
in positions $(1,n)$, $\ldots$, $(n-1,n)$ and the tuple is 
reducible again.

The proposition is proved.~~~~~$\Box$

{\bf Proof of Proposition~\ref{VmapstoW}:}

$1^0$. The mapping $\phi$ is analytic. Indeed, one can fix the contours 
defining the generators of 
$\pi _1({\bf C}P^1\backslash \{ a_1,\ldots ,a_{p+1}\} )$ which define the 
monodromy operators $M_1$, $\ldots$, $M_{p+1}$. In some neighbourhoods of the 
contours there are no poles and the dependence of the monodromy operators 
on the matrices-residua is analytic, i.e. $\phi$ is analytic.

$2^0$. Suppose that the tuples with trivial centralizers of 
matrices-residua of two different systems 
(\ref{Fuchs}) are mapped onto one and the same tuple of monodromy 
operators $M_j$. Hence, the systems are obtained from one another by a 
linear change $X\mapsto W(t)X~(*)$ of the variables $X$ where the matrix 
$W(t)$ is holomorphic and holomorphically invertible outside the poles $a_j$ 
and is at most meromorphic at $a_j$. Indeed, one has $W=X_1(X_2)^{-1}$ 
where $X_i$ are fundamental solutions of the two systems; these solutions 
and their determinants grow no faster than some power of $(t-a_j)$ when 
$t\rightarrow a_j$. Hence, $W$ (and in the same way $W^{-1}$) 
is at most meromorphic at $a_j$. Outside $a_j$ both $X_1$ and $X_2$ are 
holomorphic and holomorphically invertible, hence, $W$ as well. 

$3^0$. Under the change $(*)$ the system d$X/$d$t=A_1(t)X$ changes into 
d$X/$d$t=A_2(t)X$ where $A_2(t)=-W^{-1}($d$W/$d$t)+W^{-1}A_1(t)W$. Hence, 
one must have d$W/$d$t=A_1(t)W-WA_2(t)$. Set 
$W=W_k/(t-a_j)^k+W_{k-1}/(t-a_j)^{k-1}+\ldots$ and suppose that $k\geq 1$. 
Denote by $A_j^i$ the 
matrices-residua of the two systems, $i=1,2$. Then one has 
$-kW_k=A_j^1W_k-W_kA_j^2$, i.e. $(A_j^1+kI)W_k-W_kA_j^2=0$. The two matrices 
$A_j^1+kI$ and $A_j^2$ have no eigenvalue in common; 
indeed, one has 
$A_j^1, A_j^2\in c_j$ and $c_j$ has no eigenvalues differing by a 
non-zero integer. Hence, $W_k=0$ (see about the matrix equation $AX-XB=0$ 
in \cite{Ga}), i.e. 
$W$ has no pole at $a_j$ for any $j$ and $W$ is a constant non-degenerate 
matrix. 

$4^0$. The triviality of the centralizer of the tuple of matrices 
$A_j$ implies that $\phi$ is bijective (globally, from ${\cal V}$ onto its 
image in ${\cal W}$). An analytic bijective mapping is 
analytically invertible. This proves part 1) of the proposition. If the 
eigenvalues are generic, then every tuple of matrices 
$A_j$ from ${\cal V}$ is irreducible, hence, with 
trivial centralizer from where part 2) of the proposition follows.~~~~$\Box$

{\bf Proof of Proposition~\ref{directsum}:}

$1^0$. Prove part 1) of the proposition. Two cases are possible: 

{\em Case 1) The centralizer contains a matrix with at least two distinct 
eigenvalues.} Then the latter can be conjugated to a block-diagonal form 
with two diagonal blocks which have no eigenvalue in common. The commutation 
relations with such a matrix imply that the matrices $A_j$ or $M_j$ 
are themselves block-diagonal, i.e. they define a direct sum.

{\em Case 2) Each matrix from the centralizer is with a single 
eigenvalue.} Hence, the centralizer contains a nilpotent matrix $N$ and taking 
powers of it one can assume that $N^2=0$. One can conjugate $N$ to the 
form $\left( \begin{array}{ccc}0&0&I\\0&0&0\\0&0&0\end{array}\right)$ or 
$\left( \begin{array}{cc}0&I\\0&0\end{array}\right)$. Hence, the matrices 
$A_j$ (or $M_j$) are blocked as follows: 
$\left( \begin{array}{ccc}L&Q&R\\0&S&T\\0&0&L\end{array}\right)$ or 
$\left( \begin{array}{cc}L&R\\0&L\end{array}\right)$. The presence of two 
equal diagonal blocks implies that $A_{p+1}$ (or $M_{p+1}$) 
has a multiple eigenvalue -- a contradiction. 

Hence, one is necessarily in Case 1) and the representation defined by the 
matrices $A_j$ (or $M_j$) is a direct sum.

$2^0$. Prove part 2). If a matrix $Z$ from the centralizer is with non-trivial 
Jordan structure (i.e. has at least one Jordan block of size $>1$), then 
a suitable polynomial of $Z$ is a nilpotent matrix $N$, $N^2=0$, and like 
in case 2) one shows that $A_{p+1}$ or $M_{p+1}$ must have a multiple 
eigenvalue which is a contradiction.~~~~$\Box$

{\bf Proof of Proposition~\ref{dirsemidir}:}

$1^0$. Prove part 1)  
in the case of matrices $A_j$ (the case of matrices $M_j$ is left for the 
reader). Note that one has 
Ext$^1(P_1,P_2)=R/Q$ (which implies that 
dim~Ext$^1(P_1,P_2)=$dim$R-$dim$Q$) where 

\[ R=\{ (A_1^1X_1-X_1A_1^2,\ldots ,A_{p+1}^1X_{p+1}-X_{p+1}A_{p+1}^2)|
X_j\in M_{m_1,m_2}({\bf C}), \sum _{j=1}^{p+1}(A_j^1X_j-X_jA_j^2)=0\}\]
\[ Q=\{ (A_1^1X-XA_1^2,\ldots ,A_{p+1}^1X-XA_{p+1}^2)|
X\in M_{m_1,m_2}({\bf C})\}\]

Indeed, if $A_j=\left( \begin{array}{cc}A_j^1&0\\0&A_j^2\end{array}\right)$, 
$Y_j=\left( \begin{array}{cc}I&X_j\\0&I\end{array}\right)$, then for the 
matrix $A_j'=(Y_j)^{-1}A_jY_j$ one has 
$A_j'=\left( \begin{array}{cc}A_j^1&A_j^1X_j-X_jA_j^2\\0&A_j^2
\end{array}\right)$.  
Thus the space $R$ (resp. $Q$) is the one of tuples of right 
upper blocks of matrices $(Y_j)^{-1}A_jY_j$ whose sum is $0$ (resp. of 
tuples of matrices $Y^{-1}A_jY$, i.e. right upper blocks obtained 
as a result of simultaneous conjugation by $Y$ of the form of $Y_j$). 

On the other hand, the tangent space to the conjugacy class of the matrix 
$A_j$ (resp. of $A_j^1$ or of $A_j^2$) is 
$\{ [A_j,V_j]|V_j\in gl(m_1+m_2,{\bf C})\}$ 
(resp. the latter's restriction to the left upper or to the right 
lower block). The dimensions of these tangent spaces equal respectively 
$d_j'$, $d_j^1$, $d_j^2$. The 
restrictions to the non-diagonal blocks being of equal dimension (hence, this 
dimension equals $(d_j'-d_j^1-d_j^2)/2$), one deduces 
that dim$R=(\sum _{j=1}^{p+1}(d_j'-d_j^1-d_j^2)/2)-m_1m_2$. 

Subtracting $m_1m_2$ corresponds to the condition the sum of the matrices 
$A_j^1X_j-X_jA_j^2$ to equal $0$. This condition 
is equivalent to a system of $m_1m_2$ linear and 
linearly independent equations. Indeed, their linear dependence would 
imply the existence of a non-zero $m_2\times m_1$-matrix $S$ such that 

\[ {\rm tr}(S(\sum _{j=1}^{p+1}(A_j^1X_j-X_jA_j^2)))={\rm tr}
(\sum _{j=1}^{p+1}(SA_j^1-A_j^2S)X_j)=0\] 
identically in the entries of 
the matrices $X_j$. Hence, for $j=1\ldots ,p+1$ one has $SA_j^1-A_j^2S=0$. 
For $j=p+1$ this implies that $S=0$ because $A_{p+1}^1$ and $A_{p+1}^2$ have 
no eigenvalue in common -- a contradiction.
Obviously, dim$Q=m_1m_2$ (the map $X\mapsto A^1_{p+1}X-XA_{p+1}^2$ is 
bijective because $A_{p+1}^1$ and $A_{p+1}^2$ have 
no eigenvalue in common) which implies the formula for $\xi$.   

$2^0$. Part 2) is trivial. To deform continuously a direct 
sum into a semidirect one it suffices to replace $X_j$ by $\varepsilon X_j$, 
$\varepsilon \in ({\bf C},0)$.  

$3^0$. Prove part 3). The dimension of the variety 
${\cal V}(c_1',\ldots ,c_{p+1}')$ (resp. 
${\cal V}(c_1^i,\ldots ,c_{p+1}^i)$, $i=1,2$) where $c_j'$ (resp. $c_j^i$) 
is the conjugacy class of $A_j$ (resp. $A_j^i$) 
equals $\nu '=\sum _{j=1}^{p+1}d_j'-(m_1+m_2)^2+1$ (resp. 
$\nu ^i=\sum _{j=1}^{p+1}d_j^i-m_i^2+1$), see 
Proposition~\ref{smooth}. 

Show that if $\xi \geq 2$, then $\eta <\nu '$ where 
$\eta$ is the dimension of the variety ${\cal G}$ 
of semidirect sums of $P_1$ and 
$P_2$. One has $\eta =\eta _0+m_1m_2$ where $\eta _0$ is the dimension of 
the variety ${\cal G}_0$ 
of block upper-triangular matrices $A_j\in c_j'$ whose diagonal blocks define 
the representations $P_i$. Adding $m_1m_2$ corresponds to the possibility 
to obtain every tuple of matrices from ${\cal G}$ by conjugating 
a tuple from ${\cal G}_0$ by a matrix 
$\left( \begin{array}{cc}I&0\\X&I\end{array}\right)$ where 
$X\in M_{m_2,m_1}({\bf C})$. 

One has $\eta _0=\nu ^1+\nu ^2+$dim$R$. We let the reader check oneself that 
the inequality $\xi >1$ is equivalent to $\eta <\nu '$. The last inequality 
implies the existence of irreducible representations.~~~~$\Box$

\section{Proof of Theorem \protect\ref{basictm}\protect\label{prtm}}

{\bf Proof of part 1):}

$1^0$. Consider the case of matrices $A_j$. Multiplying by 
$b\in {\bf C}^*$ or adding to $A_j$ of scalar matrices with zero sum 
changes an (ir)reducible tuple of matrices $A_j$ 
with zero sum into such a tuple, therefore it suffices to prove the 
theorem only in the case when $c_j\subset sl(n,{\bf C})$ and there is no 
non-zero integer difference between 
two eigenvalues (for every conjugacy class $c_j$). 

Consider the mapping 
$\tau :(c_1\times \ldots \times c_p)\rightarrow sl(n,{\bf C})$, 
$\tau :(A_1,\ldots ,A_p)\mapsto -A_1-\ldots -A_p$.  
Its graph $G$ is naturally isomorphic to $(c_1\times \ldots \times c_p)$, 
hence, it is a connected smooth variety. Its subset $G_0$ consisting of 
$p$-tuples for which the matrix $A_{p+1}=-A_1-\ldots -A_p$ is with distinct 
eigenvalues and the eigenvalues of the tuple $(A_1,\ldots ,A_{p+1})$ 
are generic, is a Zariski open dense subset of $G$; hence, $G_0$ is connected. 

$2^0$. Consider the projection $\pi :G\rightarrow {\bf C}^{n-1}$ where 
${\bf C}^{n-1}$ is the space of symmetric functions of the eigenvalues of 
$A_{p+1}$. For generic 
eigenvalues the fibres are non-empty smooth varieties; a priori they can 
consist of several components which are of one and the same dimension, and 
this dimension is the same for all generic eigenvalues. 
This follows 
from Theorem~\ref{CSimpson} (and the lines that follow it) 
and Proposition~\ref{smooth}. The number of these 
components (denoted by $\chi$) is one and the same outside a proper 
algebraic subset 
$H\subset {\bf C}^{n-1}$. Hence, ${\bf C}^{n-1}\backslash H$ is connected. 

Denote by $\bar{{\bf C}}^{n-1}$ the space of eigenvalues of the matrix 
$A_{p+1}$, by $\beta$ the map ``eigenvalues'' 
$\mapsto$ ''symmetric functions of them'' 
and by $\Lambda$ the image (by $\beta$) in ${\bf C}^{n-1}$ of a hyperplane in 
$\bar{{\bf C}}^{n-1}$ defined by an 
equation $\lambda _{i_1,p+1}=\lambda _{i_2,p+1}$ (recall that by 
$\lambda _{i,k}$ we denote the 
eigenvalues of $A_k$) and by $K$ the image in ${\bf C}^{n-1}$ of 
a hyperplane in $\bar{{\bf C}}^{n-1}$ defined by some 
non-genericity relation.  Denote by 
$\tilde{\Lambda }$ (resp. $\tilde{K}$) the union of all sets  
$\Lambda$ (resp. $K$). Consider a closed oriented contour 
$\gamma \subset {\bf C}^{n-1}\backslash 
(\tilde{\Lambda }\cup \tilde{K}\cup H)$. The 
{\em monodromy} of $\gamma$ is by definition the permutation operator acting 
on the components of the fibre $\pi ^{-1}(a)$, $a\in \gamma$, when $\gamma$ 
is run once from $a$ to $a$. Up to conjugacy, this operator does not depend 
on the choice of $a\in \gamma$. 

$3^0$. We show that the monodromy of every such contour $\gamma$ is trivial. 
This together with the connectedness of $G_0$ implies that each fibre over 
${\bf C}^{n-1}\backslash (\tilde{\Lambda }\cup \tilde{K}\cup H)$ consists of 
a single component which proves the 
connectedness of the varieties ${\cal V}$ for 
such eigenvalues. 

\begin{lm}\label{H}
The fibres over 
$H\backslash (\tilde{\Lambda }\cup \tilde{K})$ are also non-empty and 
connected. 
\end{lm}

All lemmas from this section are proved in Section~\ref{prlemmas}.

Proving the triviality of the monodromy of each contour $\gamma$ is 
tantamount to showing that the monodromy of each small lace 
$\gamma '\subset {\bf C}^{n-1}$ around each 
set $\Lambda$ (or $K$) in ${\bf C}^{n-1}$ and around each component of 
$H$ of codimension $1$ in ${\bf C}^{n-1}$ is trivial. 
The lace is presumed to be 
contractible to a point $\Omega$ 
from $\tilde{\Lambda}$ (or $\tilde{K}$, or $H$). It suffices to consider 
two cases:

A) $\Omega \in \tilde{\Lambda}\backslash \tilde{K}$; we do not subtract 
$H$ because we cannot 
claim that $H$ does not contain components of $\tilde{\Lambda}$ of 
maximal dimension; it is sufficient to consider only the case when 
$A_{p+1}$ has a single double eigenvalue, its other 
eigenvalues being simple;

B) $\Omega \in (\tilde{K}\cup H)\backslash \tilde{\Lambda}$; in this case 
$A_{p+1}$ has distinct eigenvalues.

$4^0$. Consider Case B) first (in $4^0$ -- $5^0$). The eigenvalues 
of $A_{p+1}$ are defined by 
their symmetric functions only up to permutation. Fix these 
eigenvalues at 
$\Omega$, i.e. choose one of their permutations. By continuity this defines 
such a permutation for all points from ${\bf C}^{n-1}$ sufficiently close 
to $\Omega$ because the eigenvalues are distinct. In other terms, one can 
consider the covering 
$\beta :\bar{\bf C}^{n-1}\rightarrow {\bf C}^{n-1}$. There are 
$n!$ sheats that cover the neighbourhood of any point from 
${\bf C}^{n-1}\backslash \tilde{\Lambda}$. Fixing the permutation of the 
eigenvalues means fixing one of these sheats. 

For each $n$-vector $\vec{v}$ with integer components whose sum is $0$ shift 
the eigenvalues of $A_{p+1}$ by $\vec{v}$. This shifts their symmetric 
functions by some vector $\vec{v}'$ 
(recall part D) of Subsection~\ref{knownfacts}). 
Denote by $\gamma '(\vec{v}')$ 
(resp. by $c_{p+1}(\vec{v})$) the lace 
$\gamma '$ shifted by $\vec{v}'$ (resp. the semi-simple conjugacy class in 
$sl(n,{\bf C})$ whose 
eigenvalues are shifted w.r.t. the ones of $c_{p+1}$ by $\vec{v}$). Hence, 
there exists $\vec{v}=\vec{v}_0$ such that $\gamma '(\vec{v}_0')$ is 
contractible in 
${\bf C}^{n-1}\backslash (\tilde{\Lambda }\cup \tilde{K}\cup H)$. This 
follows from the algebraicity of $\tilde{\Lambda }\cup \tilde{K}\cup H$. 
Hence, the monodromy of $\gamma '(\vec{v}_0')$ is trivial. 

$5^0$. Show that the monodromies of $\gamma '$ and $\gamma '(\vec{v}_0')$ are 
the same. 
By deforming $\gamma '$ in 
${\bf C}^{n-1}\backslash (\tilde{\Lambda }\cup \tilde{K}\cup H)$ one can 
achieve the condition that for every vector $\vec{v}_0$ and for every point of 
$\gamma '(\vec{v}_0')$ the matrix $A_{p+1}$ has no eigenvalues differing by 
a non-zero integer.

Denote by ${\bf C}'$ the space of symmetric functions of the 
eigenvalues of the matrix 
$M_{p+1}=(M_1\ldots M_p)^{-1}$ where $M_j\in C_j$ for $j\leq p$. 
Hence, ${\bf C}'={\bf C}^n\cap \{ \sigma =1\}$, where  
$\sigma =\sigma _{1,p+1}\ldots \sigma _{n,p+1}$,  
and ${\bf C}'$ is a non-singular variety. Define its subvarieties 
$\tilde{\Lambda }^*$, $\tilde{K}^*$ by analogy with $\tilde{\Lambda }$, 
$\tilde{K}$ (see $2^0$) -- consider the map $\zeta :$ ``eigenvalues'' 
$\mapsto$ ``symmetric functions of them''. The set 
$\tilde{\Lambda }^*$ (resp. $\tilde{K}^*$) is the 
union of the intersections with ${\bf C}'$ of the images by $\zeta$ 
in ${\bf C}^n$ of hyperplanes $\sigma _{j_1,p+1}=\sigma _{j_2,p+1}$ 
(resp. of the varieties defined by non-genericity relations). 

There exists a proper algebraic subset $H^*\subset {\bf C}'$ such that for 
$(\sigma _{1,p+1},\ldots ,\sigma _{n,p+1})\in 
{\bf C}'\backslash (\tilde{\Lambda }^*\cup \tilde{K}^*\cup H^*)$ the sets 
${\cal W}$ are smooth varieties with one and the same 
number $\chi ^*$ of components. Notice that their dimension is the same and 
it is one and the same for all such sets of 
eigenvalues of $M_{p+1}$. This dimension 
is the same as the dimension of the components of 
${\cal V}$ for eigenvalues of $c_{p+1}$ from 
${\bf C}^{n-1}\backslash (\tilde{\Lambda }\cup \tilde{K}\cup H)$, see 
Proposition~\ref{smooth}.   

\begin{lm}\label{chi}
One has $\chi ^*=\chi$.   
\end{lm}

The mapping $\phi$ (see part E) of Subsection~\ref{knownfacts}) 
defines in a natural way images in ${\bf C}'$ of the 
contours $\gamma '$ and $\gamma '(\vec{v}_0)$. These images 
coincide; denote this common image by $\gamma ^*$. The monodromies of 
$\gamma '$ and $\gamma '(\vec{v}_0)$ are the same as the monodromy of 
$\gamma ^*$. Hence, these monodromies coincide and the monodromy 
of $\gamma '$ is trivial.

$6^0$. Consider Case A). Choose like in $4^0$ a permutation of the 
eigenvalues of $A_{p+1}$ at $\Omega$ (the permutations are $n!/2$ because 
there is one double eigenvalue). Denote by $\bar{\Omega}\in \bar{\bf C}^{n-1}$ 
the point defined by such eigenvalues. Any neighbourhood  
of $\bar{\Omega}$ 
defines a two-fold covering of some neighbourhood of $\Omega$ via $\beta$. 

The set ${\cal V}$ at $\bar{\Omega}$ is a smooth non-empty algebraic variety. 
A point from it is a tuple of matrices $A_j$ with zero sum. If in this 
tuple the matrix $A_{p+1}$ is diagonalizable, then it 
can be analytically deformed into another one from ${\cal V}$ 
with ``more generic'' JNF $J(A_{p+1})$, 
i.e. one in which $A_{p+1}$ has a single 
Jordan block of size $2$ corresponding to its double eigenvalue 
(see the description of the 
basic technical tool in part B) of Subsection~\ref{knownfacts} and 
Example~\ref{JNFnongeneric} there). So assume that $A_{p+1}$ has this ``more 
generic'' JNF. Then deform the tuple into one where $A_{p+1}$ has 
distinct eigenvalues (see Example~\ref{newEVs}). The 
parameters of the deformation can be chosen to   
belong to a neighbourhood $\omega$ of $\bar{\Omega}$ in some 
one-dimensional variety transversal at $\bar{\Omega}$ 
to the lifting in $\bar{{\bf C}}^{n-1}$ 
of $\tilde{\Lambda}$. 

Thus one can identify the components of the variety 
${\cal V}$ at $\bar{\Omega}$ with some or all components of the varieties 
${\cal V}$ at the other points from $\omega$. 
A priori the varieties 
${\cal V}$ at points from $\omega \backslash \bar{\Omega}$ have no less 
(say, $k_0$) components than the one at $\bar{\Omega}$. 

{\em The union of varieties ${\cal V}$ over 
$\omega \backslash \bar{\Omega}$ defines a trivial fibration over 
$\omega \backslash \bar{\Omega}$. This implies that the monodromies 
along laces from Case A) are 
also trivial.}   

To prove the claim use similar ideas to the ones from Case B), see $4^0$. 
Namely, there exists a shift of the eigenvalues of $A_{p+1}$ 
defined at a point from $\omega$ 
by a vector $\vec{v}$ (with integer components whose sum is $0$) so that 
the shifted neighbourhood $\omega (\vec{v})$ 
does not intersect with the lifting in 
$\bar{\bf C}^{n-1}$ of the set $\tilde{\Lambda}\cup \tilde{K}\cup H$. Hence, 
all varieties ${\cal V}$ at points from $\omega (\vec{v})$ have $k_0$ 
components. They form a trivial fibration over $\omega (\vec{v})$ (to see 
this use again the basic technical tool).

The mapping $\phi$ (see part E) of Subsection~\ref{knownfacts}) defined 
for eigenvalues from $\omega (\vec{v})$ (we denote this mapping by $\phi _1$) 
allows one to identify the components of the varieties ${\cal V}$ over 
$\omega (\vec{v})$ and the varieties 
${\cal W}$ over $\phi (\omega (\vec{v}))$. Hence, the components 
of the varieties ${\cal W}$ over $\phi (\omega (\vec{v}))$ define a trivial 
fibration over $\phi (\omega (\vec{v}))$. 

On the other hand, one can 
define the mapping $\phi$ (and at least almost everywhere in the varieties 
${\cal W}$, i.e. on the 
complements of proper analytic subsets, the mapping $\phi ^{-1}$) 
for eigenvalues from $\omega \backslash \bar{\Omega}$; denote this mapping by 
$\phi _2$. 

This means that 
the components of the varieties ${\cal V}$ define a trivial fibration 
over $\omega \backslash \bar{\Omega}$. Indeed, the varieties ${\cal V}$ over 
points from $\omega (\vec{v})$ and the respective points from $\omega$ 
(i.e. shifted by $-\vec{v}$) define the matrices-residua of 
Fuchsian systems with one and the same monodromy; hence, there exists a linear 
change $X\mapsto W(t)X$ (which is the map $\phi _2^{-1}\circ \phi _1$) 
that transforms the first systems into the second 
ones. This change defines a birational diffeomorphism between 
varieties ${\cal V}$ 
over $\omega (\vec{v})\backslash \bar{\Omega}(\vec{v})$ and 
$\omega \backslash \bar{\Omega}$ (where $\bar{\Omega}(\vec{v})$ is 
$\bar{\Omega}$ shifted by $\vec{v}$; one excludes 
$\bar{\Omega}$ because it is not clear whether the number of components of 
the variety ${\cal V}$ over $\bar{\Omega}$ is the same as for points from 
$\omega \backslash \bar{\Omega}$). The triviality of the fibration over 
$\omega (\vec{v})\backslash \bar{\Omega}(\vec{v})$ implies the one over 
$\omega \backslash \bar{\Omega}$. 

Hence, the monodromies along laces from Case A) are trivial. 
This proves part 1) of the theorem for matrices $A_j$.   

$7^0$. The connectedness of the varieties ${\cal W}$ which 
are the fibres over 
${\bf C}'\backslash (\tilde{\Lambda }^*\cup \tilde{K}^*\cup H^*)$ follows from 
Lemma~\ref{chi} (we just showed that $\chi =1$). The fibres over 
$H^*\backslash (\tilde{\Lambda }^*\cup \tilde{K}^*)$ are non-empty, see 
Theorem~\ref{CSimpson}. Their connectedness is proved by analogy with 
Lemma~\ref{H} and we leave this for the reader.
 
Part 1) of the theorem is proved. 

{\bf Proof of part 2):}

$8^0$. {\em Every tuple of matrices $A_j$ from the algebraic 
closure of ${\cal V}$ and with non-trivial centralizer can be 
continuously deformed into one from ${\cal V}$ and 
with trivial centralizer. The same is true for the varieties  
${\cal W}$.} 

The statement implies that the algebraic closure of 
${\cal V}$ belongs to its topological closure (for 
${\cal W}$ the proof of the statement is analogous and we skip it off). 
On the other hand, the topological closure belongs to the algebraic one -- 
the eigenvalues of the classes $c_j$ remain the same on the topological 
closure of ${\cal V}$ 
(by continuity); the ranks of the matrices $(A_j-\lambda _{k,j}I)^s$ 
remain the same on ${\cal V}$ and on its topological 
closure they can only drop; hence, on the topological closure of 
${\cal V}$ the matrices $A_j$ belong to the topological 
closures of the conjugacy classes $c_j$ which are also their algebraic 
closures; the condition the sum of the matrices to be $0$ is preserved (by 
continuity) on the topological closure; 
all this means that a tuple of matrices from the topological closure of 
${\cal V}$ belongs to its algebraic closure. 

To prove the statement we use the following lemmas:

\begin{lm}\label{deformclosure}
A tuple of matrices $A_j$ from the algebraic closure of 
${\cal V}$ can be continuously deformed into one from ${\cal V}$. The same 
is true for varieties ${\cal W}$.
\end{lm}

\begin{lm}\label{geq3}
Suppose that the tuples of matrices $A_j^i$, $i=1,2$, define two 
representations $P_i$ (of ranks $m_i$) with trivial centralizers where 
$\sum _{j=1}^{p+1}A_j^i=0$, $i=1,2$, and the matrices $A_{p+1}^1$, 
$A_{p+1}^2$ have each distinct eigenvalues and no eigenvalue in common. 
Then there exist representations (defined by $p+1$ matrices whose sum is $0$) 
which are semi-direct sums of $P_1$ and $P_2$ (in both possible orders) and 
with trivial centralizers at least in the following cases: 

1) $m_1\geq 3$, $m_2\geq 2$; 

2) $m_1=m_2=2$ and for at least one index $j$, $1\leq j\leq p$, 
the matrices $A_j^1$ and $A_j^2$ belong to different conjugacy classes;

3) $m_1=m_2=1$ and for at least two indices $j$, $1\leq j\leq p$, the 
matrices $A_j^1$ and $A_j^2$ are different; 

4) $m_1=2$, $m_2=1$ and for at least one index $j$, $1\leq j\leq p$, the 
matrices $A_j^i$ have no eigenvalue in common; 

5) $m_1>1$, $m_2=1$ and $r_1^1+\ldots +r_p^1>m_1$ where $r^1_j=r(A_j^1)$;

6) $m_1=m_2=2$ and for at least three indices $j$, $1\leq j\leq p$, at least 
one of the matrices $A_j^i$, $i=1,2$, is not scalar.

A similar statement holds for matrices $M_j^i$, $M_1^i\ldots M_{p+1}^i=I$. 
\end{lm}

Given a tuple of matrices from ${\cal V}$ and 
with non-trivial centralizer conjugate it to a block-diagonal form the 
diagonal blocks defining representations with trivial centralizers, see 
Proposition~\ref{directsum}. If one can 
find a couple of such representations for which Lemma~\ref{geq3} is 
applicable, then one can replace this couple by its semi-direct sum 
(which is with trivial centralizer) and thus reduce the number of 
diagonal blocks. When 
this number is $1$, then the representation is with trivial centralizer. So 
describe all cases when the lemma is not applicable:

{\em Case A)} {\em All diagonal blocks are of sizes $1$ and/or $2$.} 

{\em Case B)} {\em There is one diagonal block of size $>2$ and the rest are 
of size $1$.} 

Indeed, there can be not more than one block of size $>2$, otherwise we are 
in case 1) of the lemma, and if there is such a block, then there are no 
blocks of size $2$ for the same reason. 

\begin{rem}\label{onedim}
In both cases A) and B) if there is 
more than one diagonal block of size $1$, then for all indices $j=1,\ldots ,p$ 
but one (say, but for $j=1$) 
the restrictions of the matrices $A_j$ to all diagonal blocks of size 
$1$ are the same, otherwise case 3) of Lemma~\ref{geq3} is applicable.
\end{rem} 

In case A) for $j=1,\ldots ,p$ the conjugacy classes of the restrictions 
of the matrices $A_j$ to all blocks of size $2$ are the same (to avoid case 
2) of the lemma) and only for two of these indices (say, $1$ and $2$) are 
the restrictions of the matrices $A_j$ to the diagonal blocks of size $2$ 
non-scalar (to avoid case 6)). The one-dimensional blocks must have 
eigenvalues which for $j=1,\ldots ,p$ are eigenvalues of the 
two-dimensional blocks (to avoid case 4)) and Remark~\ref{onedim} holds. 
Hence, the matrices (if any) $A_3$, $\ldots$, $A_p$ are scalar and one can 
assume that $p=2$; each of the matrices $A_1$, $A_2$ is either 
diagonalizable, with two eigenvalues, or is with a single eigenvalue and 
with Jordan blocks of sizes $\leq 2$. Hence, for $j=1,2$ one has 
$d_j\leq n^2/2$ and as $d_3=n^2-n$, condition $(\alpha _n)$ does not hold for 
$n>2$. Thus case A) is possible only for $n=2$ and in this case either 
there is a single diagonal block and the centralizer is trivial or one 
can construct a semi-direct sum of the two one-dimensional blocks (one is 
in case 3) of the lemma, otherwise condition $(\alpha _2)$ fails).  

In case B) for $j=2,\ldots ,p$ the restrictions of $A_j$ to the 
one-dimensional blocks are equal (see Remark~\ref{onedim}). Moreover, for 
$j=1,\ldots ,p$ these restrictions are 
eigenvalues $\lambda _j$ of the restrictions $A_j^0$ of $A_j$ to the block 
of size $m_1>1$ for which one has 
rk$(A_j^0-\lambda _jI)=r(A_j^0)$. (Indeed, if not, then for the representation 
$P^0$ defined by the matrices $A_j^0$ and for some diagonal block $P^1$ 
one will have dim~Ext$^1(P^0,P^1)>r(A_1^0)+\ldots +r(A_p^0)+m_1-2m_1\geq 0$ 
and one can construct a semi-direct sum of $P^0$ and $P^1$; compare with 
case 5) of the lemma.) But then 
one has $r_1+\ldots +r_p<n$, i.e. condition $(\beta _n)$ does not hold 
for the matrices $A_j$ (the reader should go into the details oneself). 
Hence, case B) is also impossible. 
    
$9^0$. {\em If a variety ${\cal V}$ contains irreducible 
tuples, then its algebraic and topological closure is the one of 
its subset of irreducible tuples. The latter is a connected 
smooth variety. The same statements hold for varieties 
${\cal W}$.}

Fix distinct points $a_j\in {\bf C}$, $j=1,\ldots ,p+1$ which will be 
poles of Fuchsian systems. 
By multiplying the matrices 
by $\alpha \in {\bf C}^*$ one can achieve the 
condition no sum of some of the eigenvalues of the matrices $A_j$ and 
no difference $\lambda _{k_1,j}-\lambda _{k_2,j}$ to be a 
non-zero integer. This implies (see Theorem 5.1.2 from \cite{Bo1}) 
that if a tuple of matrices $A_j$ is 
irreducible, then the monodromy group of a Fuchsian system with poles $a_j$ 
is irreducible as well. 

Define the variety ${\cal W}(C_1,\ldots ,C_{p+1})$ 
after ${\cal V}(c_1,\ldots ,c_{p+1})$ like in 
Subsection~\ref{knownfacts}, part E), i.e. 
$C_j$ is the conjugacy class of the matrix  $M(A_j)$. 

Denote by $A^1$, $A^2$ 
two tuples of matrices from ${\cal V}$, the 
first irreducible and the second with trivial centralizer. In what follows 
we allow $A^2$ to be replaced in the course of the proof by a tuple 
arbitrarily close to 
it, with trivial centralizer and from the same component of ${\cal V}$ 
as $A^2$. We show 
that the component of ${\cal V}$ to which $A^2$ 
belongs, lies in the closure of the component of $A^1$. Hence, the 
same will be true without the assumption the centralizer of $A^2$ to be 
trivial, see the statement from $8^0$. Denote by 
$M^1$, $M^2$ the images of $A^1$, $A^2$ 
in ${\cal W}$ by $\phi$, see 
Subsection~\ref{knownfacts}, part E).  
  
The monodromy group defined by $M^1$ can be realized by matrices-residua from 
${\cal V}(c_1,\ldots ,c_p$, $c_{p+1}(\vec{v}))$ for all but finitely many 
vectors $\vec{v}$ (see $4^0$ and Corollary~\ref{QS}). 

\begin{lm}\label{many}
The monodromy group defined by $M^2$ or a monodromy group with trivial 
centralizer close to it and from the same component of ${\cal W}$ as $M^2$, 
can be realized by 
matrices-residua from ${\cal V}(c_1,\ldots ,c_p,c_{p+1}(\vec{v}))$ for 
infinitely many vectors $\vec{v}$. 
\end{lm}

Hence, there exists a vector $\vec{v}_0$ 
for which the eigenvalues of ${\cal V}(c_1,\ldots ,c_p,c_{p+1}(\vec{v}_0))$ 
are generic and the latter variety contains tuples $A^1_*$, $A^2_*$ 
for which the corresponding tuples of monodromy operators equal 
$M^1$, $M^2$. Denote by $w$ the linear change of variables $X\mapsto W(t)X$ 
(where $W$ is meromorphic on ${\bf C}P^1$ and holomorphically invertible 
for $t\neq a_{p+1}$) which transforms a system with matrices-residua 
from ${\cal V}_1={\cal V}(c_1,\ldots ,c_p,c_{p+1})$ into 
one with matrices-residua from 
${\cal V}_2={\cal V}(c_1,\ldots ,c_p,c_{p+1}(\vec{v}_0))$. Recall that 
${\cal V}_2$ is a smooth and connected variety (this follows from 
part 1) of the theorem which is already proved). There are proper 
algebraic subsets $U_i\subset {\cal V}_i$ where $w$ or $w^{-1}$ 
is not defined. 

One can connect $A^1_*$ and $A^2_*$ in 
${\cal V}_2$ by a contour $\gamma$ such that 
all points of $\gamma$ define Fuchsian systems with trivial 
centralizers of their monodromy groups; indeed, the tuples of 
matrices-residua for which the centralizer of the monodromy group is 
non-trivial form a proper analytic subset of 
${\cal V}_2$. Moreover, one can require 
$\gamma$ to avoid the set $U_2$. Hence, $w^{-1}(\gamma )$ is a contour 
which connects $A^1$ with $A^2$ in ${\cal V}_1$ and 
which contains no tuples with non-trivial centralizers. The 
tuple $A^1$ is irreducible; this means 
that the component of ${\cal V}_1$ to which $A^2$ 
belongs, lies in the closure of the one of $A^1$. 

If $A^2$ is irreducible as well, then the above reasoning implies that 
the components of $A^1$ and $A^2$ are parts of one and the same smooth 
variety. Hence, these components coincide.

Hence, one can connect 
$M^1$ and $M^2$ by the contour $\phi (\gamma )$ which avoids 
monodromy groups with non-trivial 
centralizers; this implies (the details are left for the reader) that if 
$M^2$ is reducible, then it 
belongs to a proper subvariety of 
${\cal W}(C_1,\ldots ,C_{p+1})$ from the closure of the set of 
irreducible tuples of ${\cal W}(C_1,\ldots ,C_{p+1})$. If $M^2$ 
is irreducible, then (like above, for $A^1$, $A^2$) it belongs to one and 
the same components of the set of irreducible tuples of 
${\cal W}(C_1,\ldots ,C_{p+1})$ as $M^1$. There remains to observe that the 
reasoning 
about $M^1$, $M^2$ (including Lemma~\ref{many}) can be performed without 
defining $M^i$ as $\phi (A^i)$. Implicitly we use here the following result 
(see \cite{ArIl}, p. 132-133): {\em if the monodromy operator $M_{p+1}$ is 
diagonalizable, 
then the monodromy group is realizable by a Fuchsian system with 
$J(A_j)=J(M_j)$ for $j=1,\ldots ,p$.}

The statement is proved. 
 
{\bf Proof of part 3):}

$10^0$. Prove part 3) of the theorem for matrices $M_j$ first. Suppose 
that for the 
conjugacy classes $C_j$ conditions $(\alpha _n)$ and $(\beta _n)$ hold 
(see Proposition~\ref{d_jr_j} and Theorem~\ref{CSimpson}) and 
that $C_{p+1}$ is with distinct eigenvalues. Suppose that the classes $c_j$ 
are defined after $C_j$ so that for $A_j\in c_j$ one has 
$\exp (2\pi iA_j)\in C_j$, see Subsection~\ref{knownfacts}, part E). 
We suppose that for every $j$ there is no non-zero integer difference 
between two eigenvalues of $c_j$. In addition, one can choose the eigenvalues 
of all classes $c_j$ to be generic (because one can shift by arbitrary 
integers whose sum is $0$ the eigenvalues of $c_{p+1}$). Hence, the variety 
${\cal V}(c_1,\ldots ,c_{p+1})$ is non-empty (because the classes $c_j$ 
also satisfy conditions $(\alpha _n)$ and $(\beta _n)$, 
see Proposition~\ref{d_jr_j}, Theorem~\ref{CSimpson} and the lines that 
follow it) and so is the variety 
${\cal W}(C_1,\ldots ,C_{p+1})$. 

Every point of ${\cal W}(C_1,\ldots ,C_{p+1})$ which is the image under the 
mapping $\phi$, see Subsection~\ref{knownfacts}, part E), is a tuple 
with trivial centralizer. Indeed, 
if the centralizer is non-trivial, then the monodromy group is a direct sum 
(see Proposition~\ref{directsum}). 
The sum $s_1$ of the eigenvalues of the matrices 
$A_j$ which are exponents relative to an invariant subspace $S_1$ must be $0$. 
Indeed, by \cite{Bo2}, Lemma 3.6, this sum is a non-positive integer; 
in the case of a direct sum, there is an invariant subspace $S_2$ such that 
$S_1\oplus S_2={\bf C}^n$; the sum $s_2$ of the exponents relative to $S_2$ is 
also a non-positive integer and $s_1+s_2=0$ (the sum of all 
eigenvalues of the matrices $A_j$ is $0$; one can assume that there are no 
eigenvalues participating both in $s_1$ and $s_2$ because 
there is no non-zero integer 
difference between any two eigenvalues of any matrix $A_j$, hence, to 
each eigenvalue $\sigma$ of $M_j$ there corresponds a single eigenvalue 
$\lambda$ of $A_j$ such that $\sigma =\exp (2\pi i\lambda )$). Hence, 
$s_1=s_2=0$. However, these equalities contradict the genericity of the 
eigenvalues of the classes $c_j$. 

Hence, the weak DSP is solvable for the classes $C_j$.

$11^0$. Prove part 3) of the theorem for matrices $A_j$. It is necessary 
to prove it 
only for non-generic eigenvalues. Given the conjugacy classes 
$c_j$ (satisfying conditions $(\alpha _n)$ and $(\beta _n)$) 
fix a vector $\vec{v}$ with integer components whose sum is $0$ 
such that the conjugacy classes $c_1$, $\ldots$, 
$c_p$, $c_{p+1}(\vec{v})$ are with generic eigenvalues (see $4^0$). Hence, 
the variety ${\cal V}(c_1,\ldots ,c_p,c_{p+1}(\vec{v}))$ is non-empty and 
smooth. 

Fix a point $A\in {\cal V}(c_1,\ldots ,c_p,c_{p+1}(\vec{v}))$ and distinct 
points $a_j\in {\bf C}$, $1\leq j\leq p+1$. Fix a sequence of vectors 
$\vec{v}_i$, $i=1,\ldots ,\mu$, $\vec{v}_1=\vec{v}$, 
$\vec{v}_{\mu}=\vec{0}$ where $\vec{v}_{i+1}=\vec{v}_i+\vec{w}_i$, $\vec{w}_i$ 
being a vector two of whose components equal $1$ and $-1$ the others 
being $0$. We want the vectors $\vec{v}_i$ to satisfy 

{\bf Condition (R).} {\em If for $i=i_0$ 
the eigenvalues of the conjugacy classes 
$c_1$, $\ldots$, $c_p$, $c_{p+1}(\vec{v}_i)$ satisfy some non-genericity 
relation, then this relation is satisfied by their eigenvalues also for 
$i>i_0$.}

Denote by $P_i$ the Procedure 
$(l_i,k_i)$ which if possible to perform would shift the eigenvalues of 
$A_{p+1}$ by $\vec{w}_i$. Hence, if all procedures $P_i$ are possible 
to perform (in the prescribed order), then they would lead from $A$ 
to a point from ${\cal V}(c_1,\ldots ,c_p,c_{p+1})$. Call Procedure $Q_i$ 
the superposition $P_i\circ \ldots \circ P_1$.

If $Q_i(A)$ is an irreducible tuple but $P_{i+1}$ cannot be 
performed because the condition $c\neq 0$ does not hold (see 
Subsection~\ref{knownfacts}, part D)), then one changes a little the 
positions of the 
poles $a_j$ and one chooses an irreducible tuple close to $Q_i(A)$ 
to which one applies $P_{i+1}$, see Proposition~\ref{Procedure}.

If $Q_i(A)$ is reducible, then conjugate it to a block upper-triangular form 
with diagonal blocks defining irreducible representations. 
Denote by $B_i$ the tuple which has the same diagonal blocks as 
$Q_i(A)$ and zeros elsewhere. Continue in the same way to perform the next 
Procedures $P_i$ only to the necessary diagonal blocks (the fact that 
it will not be necessary to perform such procedures involving different 
diagonal blocks follows from Condition (R)). In the end one 
obtains a tuple from the closure of 
${\cal V}(c_1,\ldots ,c_p,c_{p+1})$ which by the statement from $8^0$ can 
be deformed into a tuple from ${\cal V}(c_1,\ldots ,c_p,c_{p+1})$ 
with trivial centralizer.

Part 3) of the theorem is proved.~~~~$\Box$

\section{Proofs of the lemmas\protect\label{prlemmas}}

{\bf Proof of Lemma~\ref{H}:}
 
The nonemptiness of the fibres follows from Theorem~\ref{CSimpson} and the 
lines that follow it. 
Suppose that such a fibre consists of at least $2$ components. Choose 
two points -- $A^1$, $A^2$ -- 
from two different components. Assume that no sum of eigenvalues of some of 
the matrices $A_j$ is a non-zero integer; this can be achieved by multiplying 
all matrices $A_j$ by $b\in {\bf C}^*$. (Such a multiplication changes the 
conjugacy classes, but it is a diffeomorphism of fibres defined for the 
different tuples of conjugacy classe.) Hence, for fixed poles $a_j$ 
of system (\ref{Fuchs}) the monodromy groups of the two systems ($F_1$) and 
($F_2$), with 
tuples of matrices-residua $A^1$, $A^2$, are irreducible (this follows 
from Theorem 5.1.2 from \cite{Bo1}). 

By Corollary~\ref{QS}, there exist 
tuples of matrices-residua $A^1_*$, $A^2_*$ from one and the same 
variety ${\cal V}(c_1,\ldots ,c_p,c_{p+1}')$ 
(where the eigenvalues of $c_{p+1}'$ 
are obtained from the ones of $c_{p+1}$ by an admissible shift, see 
Subsection~\ref{knownfacts}, part D)) such that 
the Fuchsian systems (\ref{Fuchs}) (denoted by ($F_{1*}$), ($F_{2*}$)) with 
tuples of matrices-residua $A^1_*$, $A^2_*$ have the same monodromy 
groups as systems ($F_1$), ($F_2$). 

Denote by $w$ the linear change of variables $X\mapsto V(t)X$ bringing 
systems ($F_i$) to systems ($F_{i*}$), $i=1,2$. One can connect $A^1_*$ with 
$A^2_*$ by a contour $\gamma \subset {\cal V}(c_1,\ldots ,c_p,c_{p+1}')$ 
bypassing the 
points where $w^{-1}$ is not defined (these points form a proper subvariety 
of positive codimension). Hence, $w^{-1}(\gamma )$ connects $A^1$ with $A^2$ 
in ${\cal V}(c_1,\ldots ,c_{p+1})$ which is a contradiction.~~~~$\Box$  

{\bf Proof of Lemma~\ref{chi}:}

Fix a point (i.e. a tuple of matrices $M_j$) 
from each component of the variety 
${\cal W}(C_1,\ldots ,C_{p+1})$ where the eigenvalues of $C_{p+1}$ are from 
${\bf C}'\backslash (\tilde{\Lambda }^*\cup \tilde{K}^*\cup H^*)$. There 
exists a vector $\vec{v}_1$ such that these tuples are realized as 
tuples of monodromy operators of Fuchsian systems with 
matrices-residua $A_j\in c_j$, $j\leq p$, $A_{p+1}\in c_{p+1}(\vec{v}_1)$; 
the existence of $\vec{v}_1$ follows from Corollary~\ref{QS}. 
These tuples of matrices-residua belong to different components of 
the variety ${\cal V}(c_1,\ldots ,c_{p+1}(\vec{v}_1))$; and conversely, 
points from different (from same) components of 
${\cal V}(c_1,\ldots ,c_{p+1}(\vec{v}_1))$ 
are mapped by $\phi$ into points from different (from same) components of 
${\cal W}(C_1,\ldots ,C_{p+1})$; 
this follows from Proposition~\ref{VmapstoW}.~~~~$\Box$

{\bf Proof of Lemma~\ref{deformclosure}:}

$1^0$. We prove the lemma for matrices $A_j$, for matrices $M_j$ 
the proof is much the same. Conjugate the tuple to a block-diagonal 
form with diagonal blocks defining representations with trivial centralizers, 
see Proposition~\ref{directsum}. 
In particular, there might be a single diagonal block. 

If the given tuple is not from ${\cal V}(c_1,\ldots ,c_{p+1})$ but 
from its closure, then for some $j=j_1$ the conjugacy class  
of some of the matrices $A_{j_1}$ (we denote it by $c_{j_1}^0$) 
is not $c_{j_1}$ but belongs to its 
closure (one has $j_1<p+1$ because $A_{p+1}$ is with distinct eigenvalues). 
Given a conjugacy class from $gl(n,{\bf C})$ (which defines 
a JNF $J_1$), for any conjugacy class 
containing it in its closure and defining a JNF $J_2$ we say that $J_1$ is 
{\em subordinate} to $J_2$ (notation: $J_1\prec J_2$).  

$2^0$. In what follows we use the following fact (see \cite{He}): {\em For 
$J_1$, $J_2$ as above (i.e. $J_1\prec J_2$) one has that $J_2$ is 
obtained from $J_1$ as a 
superposition of one or several operations of the form ``a couple of Jordan 
blocks with one and the same eigenvalue and of sizes $(l,s)$, $l\geq s$, 
is replaced by a couple of Jordan blocks with the same eigenvalue and 
of sizes $(l+1,s-1)$ while the other Jordan blocks (if any) remain the 
same''}. (This fact is proved in \cite{He} only for nilpotent conjugacy 
classes which implies that it is true in the general case as well.) We say 
that the above operation is of type $(l,s)$. 

\begin{ex}
Consider the family of matrices 
$A(\varepsilon )=\left( \begin{array}{cc}J'&\varepsilon D\\0&J''
\end{array}\right)$ 
where $J'$, $J''$ are two Jordan blocks with the same eigenvalue, 
of sizes $s'$, $s''$, $s'+s''=n$, and 
$\varepsilon \in ({\bf C},0)$. If $s'\geq s''$ and $D$ has a single non-zero 
entry which is in position $(s',n)$ of $A$, then for 
$\varepsilon \neq 0$, $J(A)$ 
consists of two Jordan blocks, of sizes $s'+1,s''-1$. If $s'<s''$, then for 
$D$ having a single non-zero entry, in position $(1,s'+1)$,  
$J(A)$ consists of two Jordan blocks, of sizes $s''+1,s'-1$. 
\end{ex} 

$3^0$. Fix a sequence of operations of type $(l,s)$ leading from 
$c_{j_1}^0$ to $c_{j_1}$ (and, hence, from $J(c_{j_1}^0)$ to $J(c_{j_1})$). 
Suppose that an operation of this sequence can be performed between two 
Jordan blocks corresponding to the restriction of $A_{j_1}$ to one and the 
same diagonal block. The restrictions of the matrices $A_j$ to the 
given diagonal block defining a representation with trivial centralizer, 
one can apply the basic technical tool (see part B) of 
Subection~\ref{knownfacts}) -- one deforms $A_{j_1}$ into 
$A_{j_1}+\varepsilon G$ where for $\varepsilon \neq 0$ one has 
$A_{j_1}+\varepsilon G\in c_{j_1}$ and 
one deforms analytically the restrictions $A_j^*$ of the other matrices 
$A_j$ to the given diagonal block by conjugating them so that the sum of the 
restrictions $A_j^*$ be $0$.       

If an operation of the sequence has to be performed between two Jordan blocks, 
from two different diagonal blocks, then without loss of generality one can 
assume that these are the first two diagonal blocks. To ease the notation we 
consider only the case when there are only two diagonal blocks (in the 
general case the proof is the same). The matrices $A_j$ can be presented in 
the form $A_j=\left( \begin{array}{cc}B_j&0\\0&F_j\end{array}\right)$. 
One can deform $A_{j_1}$ into $A_{j_1}+\varepsilon G\in c_{j_1}$ 
by changing only 
the right upper block of $A_{j_1}$ (use the above example). After this one 
can conjugate the matrix $A_{p+1}$ by a matrix of the form 
$\left( \begin{array}{cc}I&\varepsilon Y\\0&I\end{array}\right)$ where $Y$ 
is chosen such that  
the sum of the matrices remains $0$. This conjugation changes the 
right upper block by $\varepsilon (B_{p+1}Y-YF_{p+1})$ 
and the possibility to choose $Y$ 
follows from the fact that $B_{p+1}$ and $F_{p+1}$ have no eigenvalue in 
common and, hence, 
the linear operator $Y\mapsto B_{p+1}Y-YF_{p+1}$ is bijective. 

Hence, after finitely many operations of the form $(l,s)$ one obtains a 
tuple of matrices from ${\cal V}(c_1,\ldots ,c_{p+1})$.~~~~$\Box$
  
{\bf Proof of Lemma~\ref{geq3}:}

$1^0$. One has to check that in all 6 cases one has dim~Ext$^1(P_1,P_2)>0$. 
We do this in detail only in case 1), the most difficult one. 
Recall that one has Ext$^1(P_1,P_2)=R/Q$, hence,  
dim~Ext$^1(P_1,P_2)=$dim$R-$dim$Q$ where the spaces $R$ and $Q$ were defined 
in the proof of Proposition~\ref{dirsemidir}.  

Define the linear operator 
$\xi _j:M_{m_1,m_2}({\bf C})\rightarrow M_{m_1,m_2}({\bf C})$ by 
$\xi _j:(.)\mapsto A_j^1(.)-(.)A_j^2$. Set 
$R_j=\{ A_j^1X_j-X_jA_j^2|X_j\in M_{m_1,m_2}({\bf C})\}$. 
One has dim$R_{p+1}=m_1m_2$ because 
$A_{p+1}^i$, $i=1,2,$ have no eigenvalue in common.

Denote by $r_j^i$, $d_j^i$ the quantities 
$r_j$, $d_j$ computed for the 
matrices $A_j^i$, see the Introduction. 

$2^0$. {\em For $j\leq p$ 
one has dim$R_j\geq r_j^1m_2$. 
This is proved in $5^0$.}
 
Hence, 

\begin{equation}\label{Ext} 
{\rm dim~Ext}^1(P_1,P_2)=\sum _{j=1}^{p+1}{\rm dim}R_j-m_1m_2-{\rm dim}Q
\geq m_1m_2+(r_1^1+\ldots +r_p^1)m_2-m_1m_2-m_1m_2\geq 0
\end{equation}  
because $r^1_1+\ldots +r_p^1\geq m_1$, see Proposition~\ref{trivcentr}. 
The first term $(-m_1m_2)$ corresponds to the condition 
$\sum _{j=1}^{p+1}(A_j^1X_j-X_jA_j^2)=0$, the second equals dim$Q$ (for no 
matrix $X\in M_{m_1,m_2}({\bf C})$ does one have $A_{p+1}^1X-XA_{p+1}^2=0$ 
because $A_{p+1}^i$ have no eigenvalue in common).

$3^0$. Consider case 1) of the lemma. 
Try to understand when there might be equality in (\ref{Ext}). 
Suppose first that $A_j^i$ are diagonalizable. Hence, for $j\leq p$ at least 
two matrices $A_j^2$ (say, for $j=1,2$) are non-scalar, i.e. have at least 
two eigenvalues each, otherwise condition $(\alpha _{m_2})$ fails. 
The equality dim$R_j=r_j^1m_2$ is possible only if each 
eigenvalue of $A_j^2$ is eigenvalue of $A_j^1$ of multiplicity $m_1-r_j^1$. 
Hence, for $j=1,2$, $A_j^1$ has at least two eigenvalues of maximal 
multiplicity which is $\leq m_1/2$. 

In order to have $r_1^1+\ldots +r_p^1=m_1$, for $3\leq j\leq p$ 
the matrices $A_j^1$ must 
be scalar and the matrices $A_1^1$, $A_2^1$ must have each two eigenvalues, 
of multiplicity $m_1/2$ (i.e. $m_1$ must be even, $m_1\geq 4$). But in this 
case condition $(\alpha _{m_1})$ fails -- one has $d_{p+1}^1=m_1^2-m_1$, 
$d_1^1=d_2^1=m_1^2/2$, $d_j=0$ for $3\leq j\leq p$. Hence, for some 
$3\leq j\leq p$ there is a 
non-scalar matrix $A_j^1$ which means that $r_1^1+\ldots +r_p^1>m_1$. 

In the general case (when $A_j^i$ are not necessarily diagonal) use 
part {\bf C)} of Subsection~\ref{knownfacts}. The dimensions of the images of 
the operators $\xi _j$ are the same when they are defined after a matrix 
$A=\left( \begin{array}{cc}A_j^1&0\\0&A_j^2\end{array}\right)$ or after 
a matrix $B=\left( \begin{array}{cc}B_j^1&0\\0&B_j^2\end{array}\right)$ where 
the JNF $J(B)$ is the diagonal JNF corresponding to $J(A)$ and 
$J(B_j^i)$ corresponds to $J(A_j^i)$. Indeed, this dimension equals 
$(d_j-d_j^1-d_j^2)/2$ and the quantities $d_j$, $d_j^i$ are the same for 
both matrices $A$ and $B$.

This means that in case 1) there is always a strict inequality in (\ref{Ext}).
  
$4^0$. In cases 2), 3), 4), 5) and 6) one checks directly that 
dim~Ext$^1(P_1,P_2)>0$; in case 2) one should 
notice that for at least three indices $j$ the matrices $A_j^i$ must be 
non-scalar (for $i=1,2$), otherwise condition $(\alpha _{m_i})$ fails. We 
leave the details for the reader.

$5^0$. Without loss of generality assume that 
$A_j^2$ is in JNF. Then $\xi _j$ splits into a direct sum of operators 
acting on $M_{m_1,l_{\nu}}({\bf C})$ where $l_{\nu}$ is the multiplicity of an 
eigenvalue of $A_j^2$. So one can consider the case when there is a single 
eigenvalue $\lambda$ of $A_j^2$. Hence, 
$\xi _j:(.)\mapsto (A_j^1-\lambda I)(.)-(.)N$ where $N$ is the nilpotent 
part of $A_j^2$. A change of variables $(.)\mapsto (.)L$ with a suitable 
diagonal matrix $L$ brings the operator to the form 
$\xi _j:(.)\mapsto (A_j^1-\lambda I)(.)-(.)\varepsilon N$, 
$\varepsilon \in {\bf C}^*$. Hence, when 
$\varepsilon$ is small, $\xi _j$ is a perturbation of the operator 
$\xi _j^0:(.)\mapsto (A_j^1-\lambda I)(.)$ the dimension of whose image is 
$\geq r_j^1m_2$. Hence, the same is true for $\xi _j$ as well, i.e. 
dim$R_j\geq r_j^1m_2$.~~~~$\Box$  

{\bf Proof of Lemma~\ref{many}:}

$1^0$. Consider the Fuchsian system d$X/$d$t=(\sum _{j=1}^{p+1}A_j/(t-a_j))X$. 
Suppose that its matrices-residua are block upper-triangular: 
$A_j=\left( \begin{array}{cc}B_j&D_j\\0&R_j\end{array}\right)$. Suppose also 
that the tuples of matrices $B_j$ and $R_j$ are irreducible, that 
the only non-genericity relation modulo ${\bf Z}$ 
which holds for the eigenvalues of the 
matrices $A_j$ is the sum of the eigenvalues of the matrices $B_j$ 
(resp. $R_j$) to be $0$, and that  
$A_{p+1}=$diag$(\lambda _{1,p+1},\ldots ,\lambda _{n,p+1})$. Present 
the system by its Laurent series expansion at the pole $a_{p+1}$: 
d$X/$d$t=(A_{p+1}/(t-a_{p+1})+F+o(1))X~(*)$ where 
$F=\sum _{j=1}^pA_j/(a_{p+1}-a_j)$. 

Suppose that the block $D$ of the matrix $F$ has a non-zero entry $\beta$ 
in position $(l,k)$. In what follows the reasoning is close to the one when 
Procedure $(l,k)$ is described in part D) of Subsection~\ref{knownfacts}, 
but the monodromy group is not irreducible here. The change of variables 
$X\mapsto (I+W/(t-a_{p+1}))X$ where the only non-zero entry of $W$ is 
$W_{k,l}=(\lambda _{k,p+1}-\lambda _{l,p+1}+1)/\beta$ 
transforms the system into 
another Fuchsian system the eigenvalues of whose residuum at $a_{p+1}$ 
have changed as follows: $\lambda _{k,p+1}\mapsto \lambda _{k,p+1}+1$, 
$\lambda _{l,p+1}\mapsto \lambda _{l,p+1}-1$, the other eigenvalues 
remain the same. 

Hence, after the change of variables the new tuple of 
matrices-residua is irreducible -- the change has destroyed the only 
non-genericity relation (but, of course, has not destroyed it 
modulo ${\bf Z}$).

$2^0$. If no entry of the restriction of the matrix $F$ to the block $D$ 
is non-zero, then change a little the positions of the poles of the system so 
that one of the entries become non-zero. (This is possible to do because 
otherwise the entries of all blocks $D_j$ must be $0$, the monodromy 
group will be a direct sum and its centralizer will not be trivial.) 
The centralizer of the monodromy group remains trivial under small 
changes of the positions of the poles and the monodromy 
group remains within the same component of ${\cal W}(C_1,\ldots ,C_{p+1})$. 
One can change the poles $a_j$ so that at least one entry of the block $D$ 
becomes non-zero. After this one performs the change of variables from $1^0$. 
 

$3^0$. If the tuple of matrices $A_j$ has a more complicated block 
upper-triangular form, then one can destroy in the same way one by one all 
non-genericity relations satisfied by the eigenvalues of the matrices-residua. 
One combines changes of the positions of the poles with linear changes of the 
variables. Having destroyed all such relations, one can continue to 
construct similar 
linear changes which do not restore any of the relations. (The fact that 
the tuple of matrices-residua is no longer reducible does not 
affect the reasoning; we leave the details for the reader.) Thus one can 
construct infinitely many Fuchsian systems with $A_j\in c_j$ for $j\leq p$, 
with monodromy as required  
and with different sets of eigenvalues of $A_{p+1}$. 
The lemma is proved.~~~~$\Box$

Author's address: Universit\'e de Nice -- Sophia Antipolis, Laboratoire de 
Math\'ematiques, Parc Valrose, 
06108 Nice, Cedex 2, France; e-mail: kostov@math.unice.fr 
\end{document}